\documentclass[12pt]{amsart}

\usepackage{amsmath,amsfonts,amsthm,amsopn,amssymb}
\usepackage{verbatim,wasysym,mathrsfs}
\usepackage[left=2.6cm,right=2.6cm,top=3.1cm,bottom=3.1cm]{geometry}
\usepackage{cite,marginnote}
\usepackage{color,enumitem,graphicx}
\usepackage[colorlinks=true,urlcolor=blue, citecolor=red,linkcolor=blue,
linktocpage,pdfpagelabels, bookmarksnumbered,bookmarksopen]{hyperref}
\usepackage[english]{babel}
\usepackage[T1]{fontenc}
\usepackage{lmodern}
\usepackage[hyperpageref]{backref}

\numberwithin{equation}{section}

\makeindex
\newtheorem{theorem}{Theorem}[section]
\newcounter{maintheorem}

\newtheorem{mainth}[maintheorem]{Theorem}
\newtheorem{proposition}[theorem]{Proposition}
\newtheorem{lemma}[theorem]{Lemma}
\newtheorem{remark}[theorem]{Remark}

\newtheorem{corollary}[theorem]{Corollary}
\newtheorem{definition}{Definition}[section]

\def\e{{\rm e}}

\newcommand{\R}{\mathbb R}

\def\N{\mathbb N}

\newcommand{\tv}{\widetilde v}
\newcommand{\dd}{{\,\rm d}}
\newcommand{\bw}{\bar w}
\newcommand{\Wrad}{W_{\text{rad}}}
\newcommand{\les}{\lesssim}
\newcommand{\ges}{\gtrsim}
\newcommand{\intN}{\int_{\R^N}}
\newcommand{\supp}{{\rm supp}\,}

\newcommand{\loc}{\operatorname{\rm loc}}

\title[Nonlocal Schr\"odinger-Poisson systems in $\R^N$]{Nonlocal Schr\"odinger-Poisson systems in $\R^N$: \\the fractional Sobolev limiting case}

\author[D. Cassani]{Daniele Cassani}
\author[Z. Liu]{Zhisu Liu}
\author[G. Romani]{Giulio Romani$^*$}

\address[Daniele Cassani and Giulio Romani]{
	\newline\indent Dipartimento di Scienza e Alta Tecnologia,
	\newline\indent Universit\`{a} degli Studi dell'Insubria
	\newline\indent and
	\newline\indent RISM-Riemann International School of Mathematics
	\newline\indent Villa Toeplitz, Via G.B. Vico, 46 - 21100 Varese, Italy}
	\email{\href{mailto:daniele.cassani@uninsubria.it}{daniele.cassani@uninsubria.it}}
\email{\href{mailto:giulio.romani@uninsubria.it}{giulio.romani@uninsubria.it}}

\address[Zhisu Liu]{\newline\indent Center for Mathematical Sciences/School of Mathematics and Physics,
\newline\indent
China University of Geosciences,
\newline\indent
Wuhan, Hubei, 430074, PR China}
\email{\href{mailto:liuzhisu@cug.edu.cn}{liuzhisu@cug.edu.cn}}
\dedicatory{Dedicated to Enzo Mitidieri on the occasion of his 70$^\text{th}$ birthday}
\thanks{$^*$Corresponding author: Giulio Romani ({\tt giulio.romani@uninsubria.it.})}

\thanks{Zhisu Liu is supported by the National Natural Science Foundation of China (No.~12226331), and the Fundamental Research Funds for the Central Universities, China University of Geosciences (Wuhan, Grant number: CUG2106211; CUGST2), and Guangdong Basic and Applied Basic Research Foundation (No. 2023A1515011679). Giulio Romani is partially supported by INdAM-GNAMPA Project 2023 titled \textit{Interplay between parabolic and elliptic PDEs} (Gran number: CUP E53C2200l93000l).}

\date{\today}

\begin{document}

\begin{abstract} We study the existence of positive solutions for nonlocal systems in gradient form and set in the whole $\R^N$. A quasilinear fractional Schr\"odinger equation, where the leading operator is the $\frac Ns$-fractional Laplacian, is coupled with a higher-order and possibly fractional Poisson equation. For both operators the dimension $N\geq 2$ corresponds to the limiting case of the Sobolev embedding, hence we consider nonlinearities with exponential growth. Since standard variational tools cannot be applied due to the sign changing logarithmic Riesz kernel of the Poisson equation, we employ a variational approximating procedure for an auxiliary Choquard equation, where the Riesz kernel is uniformly approximated by polynomial kernels. Qualitative properties of solutions such as symmetry, regularity and decay are also established. Our results extend and complete the analysis carried out in the planar case in \cite{CLR}.
\end{abstract}

\maketitle

\section{Introduction}\label{sec1}

\noindent We investigate existence of solutions for the nonlocal Schr\"odinger-Poisson system in the whole space given by
\begin{equation*}\tag{SP$_s$}\label{SPs}
	\begin{cases}
		(-\Delta)_\frac Ns^su+V(x)|u|^{\frac Ns-2}u=\phi\,f(u)\\
		&\mbox{ in }\ \,\R^N,\\
		(-\Delta)^{\frac N2}\phi=F(u)
	\end{cases}
\end{equation*}
where $s\in(0,1)$, $V$ is a positive and bounded potential, $f:\R\to\R$ is a positive nonlinearity, and $F$ its primitive vanishing at zero. The nonlocal operator in the first Schr\"odinger equation is the $s-$fractional $p-$Laplacian with $p=\frac Ns$, while the operator in the second Poisson equation is the $\frac N2$-power of the Laplacian, and therefore possibly of higher-order and fractional, depending on the dimension $N\geq 2$. As we are going to see, the functional setting is critical with respect to the Sobolev embeddings, and this on the one hand yields to consider nonlinearities with exponential growth. On the other hand, the Riesz kernel of the Poisson equation is logarithmic, hence sign-changing and unbounded from below and above. We are interested in understanding the interplay between these features, together with the nonlocal and higher-order character of the operators involved.
\vskip0.2truecm
\noindent Schr\"odinger-Poisson systems emerge in several fields of Physics, such as in Hartree models for crystals, astrophysics, electromagnetism, and quantum mechanics. We refer the interested reader to \cite{BF,LRZ} for the background.

\noindent In the local semilinear case, when $N\geq3$, namely
\begin{equation*}\tag{SP}\label{SP}
	\begin{cases}
		-\Delta u+V(x)u=\phi\,f(u)\\
		&\mbox{in}\ \,\R^N,\\
		-\Delta\phi=F(u)
	\end{cases}
\end{equation*}
a well-known strategy to find solutions is to solve the Poisson equation by means of the positive polynomial Riesz kernel and substitute $\phi=\phi(u):=C_N|x|^{2-N}\ast F(u)$ into the first equation of the system. In this way, one obtains a Schr\"odinger equation with a convolution nonlinear term, the so-called Choquard type equation, which has the advantage to exhibit a variational structure and hence can be solved in the natural space $H^1(\R^N)$, see e.g. \cite{MV} and references therein. In the case $N=2$ this technique cannot be directly employed: indeed, the Riesz kernel of the Laplace operator is $\frac1{2\pi}\log\frac1{|\cdot|}$, and the nonlocal term in the Choquard equation
\begin{equation*}\tag{Ch}\label{Ch}
	-\Delta u+V(x)u=\frac1{2\pi}\left(\log\frac1{|\cdot|}\ast F(u)\right)f(u)\quad\ \mbox{in}\ \R^2
\end{equation*}
is no longer well-defined in the natural Hilbert space $H^1(\R^2)$. Let us also point out a more theoretical as well as delicate aspect: if from one side it is natural to consider the convolution with the Riesz kernel as solution of the Poisson equation, on the other side this choice is somehow arbitrary as one cannot distinguish solutions which differ by additive constants in some reasonable function space setting; see the discussion carried out in \cite{BCT,BCS,CLR}.

\noindent Following a first approach of Stubbe in \cite{Stubbe}, Cingolani and Weth managed to bypass this problem for $f(u)=u$ (strictly speaking the Choquard case) in \cite{CW} by restricting to the constrained space $\{u\in H^1(\R^2)\,|\,\int_{\R^2}\log(1+|x|)u^2\dd x<+\infty\}$, in which the functional associated to the Choquard equation turns out to be well-defined; see also further extensions in \cite{DW,CW2}.
A different approach, suitable for tuning the phenomenon of logarithmic growth of the Riesz kernel and the maximal exponential growth for the nonlinearity, was proposed by Cassani and Tarsi in \cite{CT}: in a log-weighted space built on $H^1(\R^2)$ they proved a Poho\v zaev-Trudinger inequality, which yields the well-posedness of the functional associated to \eqref{Ch} for subcritical or critical  nonlinearities $f$ in the sense of Trudinger-Moser. This approach has been then extended in \cite{BCT} to quasilinear Schr\"odinger-Poisson systems in $\R^N$, namely \eqref{SPs} with $s=1$, and in \cite{CMR} for Schr\"odinger-Poisson systems with potential and weight functions decaying to $0$ at infinity. We also refer to \cite{ACTY,AFS} for related results on Choquard equations with exponential nonlinearities and polynomial kernels, and to \cite{AF,dACS,BM2,BM3} for exponential nonlinearities and logarithmic kernels different from \eqref{SPs}. 

\noindent In all the aforementioned works the function space has to be adapted in order to have the energy functional well-defined. Recently, this difficulty was overcome in \cite{LRTZ}, see also \cite{CDL} and the extension to the zero-mass case in \cite{R}, by means of an approximating procedure. The idea is to replace the logarithmic kernel with a more suitable power-like kernel, by observing that
\begin{equation}\label{key_convergence}
	\log\frac1t=\lim_{\mu\to0^+}\frac{t^{-\mu}-1}\mu\ .
\end{equation}
\noindent Hence, one is lead to find critical points of auxiliary approximating equations with power-like convolution terms, and then retrieve a solution of the original logarithmic Choquard equation by means of a limit procedure. The advantage is that one can deal with the approximating functionals in the standard space $H^1(\R^2)$ by means of the Hardy-Littlewood-Sobolev inequality. The price to pay is the need of estimates for the family of critical points, which have to be uniform with respect to the parameter $\mu$.
\vskip0.2truecm
\noindent Short-long range interactions in the physics model yield the appearance of nonlocal operators which, from the functional analysis point of view, turn out to be well defined within fractional Sobolev-Slobodeckij spaces $W^{s,p}(\R^N)$, see for instance \cite{CVW} for more references on this topic. As in the Sobolev local case, the embedding $W^{s,p}(\R^N)\hookrightarrow L^{\frac{Np}{N-sp}}(\R^N)$ is continuous provided $N>sp$. In the borderline case $N=ps$, still one has that the maximal degree of summability of a function with membership in $W^{s,N/s}(\R^N)$ is of exponential type, as obtained by Parini and Ruf \cite{PR} for bounded domains and by Zhang \cite{Z} for the whole space, see Theorem \ref{MT_ineqPR} below. This entices into investigating the system \eqref{SPs}, where $f$ has the maximal exponential growth at infinity. To the best of our knowledge, the only results for $p$-fractional Choquard equations in $\R^N$ with exponential nonlinearities, which however do not derive from \eqref{SPs}, are obtained in \cite{CDB,BM1,YRCW,YRTZ,ChWa}. The study of \eqref{SPs} for the planar case $N=2$ has been carried out in \cite{CLR}. Since it seems difficult to obtain in this setting some log-weighted Poho\v zaev-Trudinger inequality in the spirit of \cite{CT}, here we follow the approximating method of \cite{LRTZ,CDL}. However, the effort is not only technical because of the higher dimension, as we have to deal with the basic and still unknown fact, whether the energy level reached by Moser sequences (usually equivalent to instantons or Aubin-Talenti functions in the power case) matches the sharp exponent in the fractional Trudinger-Moser inequality of \cite{PR,Z}, and this prevents to exploit concentrating sequences to gain variational compactness. The restriction to the planar case somehow simplifies technicalities and moreover, in this case the Poisson equation is semilinear and of second-order. In \cite{CLR}, we also proved that from a weak solution of the Choquard equation we can retrieve a distributional solution of the system and this is somehow neglected in the literature. See also \cite{DM,CDM,MP} for related problems and non existence results. 

\medskip

\noindent In this paper we aim at completing the analysis carried out in \cite{CLR} by considering the general system \eqref{SPs} in any dimension $N\geq2$. Note that now the Poisson equation is of higher-order and nonlocal when $N$ is odd. Before giving the precise statement of our results, let us recall some basic facts and make precise the definition of solutions we are going to use.

\noindent The $(s,\tfrac Ns)-$fractional Laplace operator pointwisely acts as
$$
(-\Delta)^s_{\frac Ns}u(x):=2\,PV\!\int_{\R^N}\frac{|u(x)-u(y)|^{\frac Ns-2}(u(x)-u(y))}{|x-y|^{2N}}\dd y,\quad x\in\R^N,
$$
where \textit{PV} stands for the Cauchy Principal Value. This is well-defined for all $x\in\R^N$ for functions in $C^{1,1}_{loc}(\R^N)$ which enjoy suitable integrability conditions at infinity, see \cite[Lemma 5.2]{CL}. One can also consider such operator in a weak form in the corresponding fractional Sobolev-Slobodeckij space
$$W^{s,\frac Ns}(\R^N):=\big\{u\in L^{\frac Ns}(\R^N): [u]_{s,\frac Ns}<+\infty\big\},$$
where $[u]_{s,\frac Ns}$ denotes the Gagliardo seminorm
$$[u]_{s,\frac Ns}=\bigg(\int_{\R^N}\int_{\R^N}\frac{|u(x)-u(y)|^{\frac Ns}}{|x-y|^{2N}}\dd x\dd y\bigg)^{\frac sN},$$
which is a uniformly convex Banach space with norm
$$\|u\|:=\bigg(\|u\|_{\frac Ns}^{\frac Ns}+[u]_{s,\frac Ns}^{\frac Ns}\bigg)^{\frac sN}.$$
Note that, if the potential $V:\R\to\R$ is measurable and satisfies
\begin{itemize}
	\item[($V$)] There exist constants $\underline{V_0},\overline{V_0}>0$ such that
	$$0<\underline{V_0}\leq V(x)\leq\overline{V_0}<+\infty\quad\ \mbox{for all}\ \,\,x\in\R^N,$$
\end{itemize}
then
$$\|u\|_V:=\bigg(\|V(\cdot)^{\frac sN}u\|_{\frac Ns}^{\frac Ns}+[u]_{s,\frac Ns}^{\frac Ns}\bigg)^{\frac sN}$$
is an equivalent norm on $W^{s,\frac Ns}(\R^N)$. In this functional setting we can define a weak notion of solution for \eqref{SPs}. For $\gamma>0$ the weighted Lebesgue space $L_\gamma(\R^N)$ is defined as
$$L_\gamma(\R^N):=\Big\{u\in L^1_{\loc}(\R^N)\,\Big|\,\int_{\R^N}\frac{|u(x)|}{1+|x|^{N+2\gamma}}\dd x<+\infty\Big\}\,.$$
As usual, we denote by $\mathcal{S}$ the Schwartz space of rapidly decreasing functions and by $\mathcal{S}'$ the dual space of tempered distributions.

\begin{definition}\label{sol_Poisson}
	For $\mathfrak f\in\mathcal S'(\R^N)$ we say that a function $\phi\in L_1(\R^N)$ is a solution of the linear Poisson equation $(-\Delta)^\frac N2\phi=\mathfrak f$ in $\R^N$ if
	$$\intN\phi(-\Delta)^\frac N2\varphi=\langle\mathfrak f,\varphi\rangle\qquad\mbox{for all}\ \,\varphi\in\mathcal S(\R^N)\,.$$
\end{definition}

\begin{definition}[Solution of \eqref{SPs}]\label{sol_SP}
	We say that $(u,\phi)$ is a weak solution of the Schr\"odinger-Poisson system \eqref{SPs} if
	\begin{equation*}
		\begin{split}
			\int_{\R^N}\int_{\R^N}&\frac{|u(x)-u(y)|^{\frac Ns-2}(u(x)-u(y))(\varphi(x)-\varphi(y))}{|x-y|^{2N}}\dd x\dd y\\
			&\qquad\qquad\qquad\quad\qquad\qquad+\int_{\R^N}V(x)|u|^{\frac Ns-2}u\varphi\dd x=\int_{\R^N}\phi f(u)\varphi\dd x
		\end{split}
	\end{equation*}
	for all $\varphi\in W^{s,\frac Ns}(\R^N)$, and $\phi$ solves $(-\Delta)^\frac N2\phi=F(u)$ in $\R^N$ in the sense of Definition \ref{sol_Poisson}.
\end{definition}

\noindent In order to find a solution of \eqref{SPs} in the sense of Definition \ref{sol_SP}, we consider the logarithmic Choquard equation
\begin{equation}\label{Chs}\tag{Ch$_s$}
	(-\Delta)^s_{\frac Ns} u+V(x)|u|^{\frac Ns-2}u=C_N\left(\log\frac1{|\cdot|}\ast F(u)\right)f(u)\quad\text{in}\ \,\R^N\,,
\end{equation}
where the Riesz kernel of the Poisson equation in \eqref{SPs}
\begin{equation}\label{Riesz_log}
	I_N(x):=C_N\log\frac1{|x|}\qquad\mbox{with}\quad C_N^{-1}:= 2^{N-1}\pi^{\frac N2}\Gamma\left(\tfrac N2\right),
\end{equation}
has been substituted in the Schr\"odinger equation in \eqref{SPs}.
\begin{definition}[Solution of \eqref{Chs}]\label{sol_Choquard_log}
	We say that $u\in W^{s,\frac Ns}(\R^N)$ is a \textit{weak solution of} \eqref{Chs} if
	\begin{equation}\label{sol_Choquard_log_test}
		\begin{split}
			\int_{\R^N}\int_{\R^N}&\frac{|u(x)-u(y)|^{\frac Ns-2}(u(x)-u(y))(\varphi(x)-\varphi(y))}{|x-y|^{2N}}\dd x\dd y+\int_{\R^N}V(x)|u|^{\frac Ns-2}u\varphi\dd x\\
			&=\int_{\R^N}\left(\int_{\R^N}\log\frac1{|x-y|}F(u(y))\dd y\right)f(u(x))\varphi(x)\dd x
		\end{split}
	\end{equation}
	for all $\varphi\in W^{s,\frac Ns}(\R^N)$.
\end{definition}

\noindent A key ingredient in our assumptions on $f$ will be the following fractional version of the Trudinger-Moser inequality obtained in \cite{PR,Z}. Let
$$\Phi_{N,s}(t):=\e^t-\sum_{j=0}^{j_{\frac Ns}-2}\frac{t^j}{j!}\,,$$
for $t\geq 0$, where $j_{\frac Ns}:=\min\{j\in\N:\,j\geq\frac Ns\}$.
\begin{mainth}[Theorem 1.3 \cite{Z}]\label{MT_ineqPR}
	Let $s\in(0,1)$, then for all $\alpha>0$ one has
	\begin{equation}\label{Trudinger}
		\int_{\R^N}\Phi_{N,s}(\alpha|u|^{\frac N{N-s}})\dd x<+\infty\,.
	\end{equation}
	Moreover, for all $\lambda>0$, set
	$$\alpha_*:=\sup\bigg\{\alpha:\,\sup_{u\in W^{s,\frac Ns}(\R^N),\,[u]_{s,\frac Ns}+\lambda\|u\|_{\frac Ns}\leq1}\int_{\R^N}\Phi_{N,s}(\alpha|u|^{\frac N{N-s}})\dd x<+\infty\bigg\}\,,$$
	one has $\alpha_*\leq\alpha_{s,N}^*\,$,
	where
	$$\alpha_{s,N}^*:=N\left(\frac{2(N\omega_N)^2\Gamma(1+\tfrac Ns)}{N!}\sum_{k=0}^{+\infty}\frac{(N-1+k)!}{k!(N+2k)^{N/s}}\right)^{\frac s{N-s}}.$$
\end{mainth}

\noindent As remarked in \cite{PR,Z}, the result of Theorem \ref{MT_ineqPR} is not sharp in the sense of Moser, since $\alpha_{s,N}^*$, which is the level reached by the concentrating fractional version of Moser sequences, is just an upper bound for the sharp exponent $\alpha_*$. Obtaining the precise value of $\alpha_*$ is still a challenging open problem. This missing information is responsible for some extra issues we have to deal with in establishing fine estimates for the mountain-pass level, see the discussion which introduces Lemma \ref{Lemma_bar w} below.
\vskip0.2truecm
\paragraph{\textbf{Assumptions.}}
According to Theorem \ref{MT_ineqPR}, the maximal growth for $f$ is exponential. We will entail this in the following set of assumptions:
\begin{itemize}	\item[($f_1$)] $f\in C^1(\R)$, $f\geq 0$, $f(t)=0$ for $t\leq0$, and $f(t)=o(t^{\frac Ns-1})$ as $t\to0^+$;
	\item[($f_2$)] there exist constants $b_1,b_2>0$
	such that for any $t>0$,
	$$0<f(t)\leq b_1+b_2\Phi_{N,s}(\alpha_*|t|^{\frac N{N-s}})\,;$$
	\item[($f_3$)] there exists $\tau\in\left(\left(1-\tfrac2N\right)s,s\right)$ such that
	$$1-s+\tau\leq\frac{F(t)f'(t)}{f^2(t)}<1+\mu_N(s,\tau)\quad\mbox{for any}\ \,t>0,$$
	where $\mu_N(s,\tau)$  is explicitly given in \eqref{def_mu};
	\item[($f_4$)] $\lim\limits_{t\to+\infty}\frac{F(t)f'(t)}{f^2(t)}=1$ or equivalently
	$\lim\limits_{t\to+\infty}\frac{\rm d}{{\rm d}t}\frac{F(t)}{f(t)}=0\,$;
	\item[($f_5$)] there exists $\beta_0>1$ depending on $s$ such that
	\begin{equation*}
		f(t)F(t)\geq\beta\,t^\frac Ns\quad\ \mbox{for all}\ \ t>T_N(s)
	\end{equation*}
	for some $\beta>\beta_0$. The values of $\beta_0$ and $T_N(s)$ are explicitly given in Lemma \ref{Lem:MP_level}.
\end{itemize}
The set of assumptions on $f$ coincides for $N=2$ with the one considered in \cite{CLR}, to which we refer for further details and comments, see also Section \ref{Sec_2}. Let us just point out that our nonlinearities can have critical or subcritical growth in the sense of Theorem \ref{MT_ineqPR}, but at least exponential, as prescribed by ($f_3$)-($f_4$). Finally, ($f_5$) is a condition in the interval $(T_N(s),+\infty)$, which requires a sufficiently fast growth for middle-range values of $t$. Indeed, at infinity, it is automatically satisfied.
\vskip0.2truecm
\paragraph{\textbf{Main results.}} We are now in a position to state our main results. First, we find a positive weak solution in the sense of Definition \ref{sol_Choquard_log} to the Choquard equation \eqref{Chs}.
\begin{theorem}[Existence for \eqref{Chs}]\label{Thm_Chs}
	Suppose $V$ is radially symmetric and satisfies (V), and that ($f_1$)-($f_5$) are fulfilled. Then, \eqref{Chs} possesses a positive radially symmetric weak solution $u\in W^{s,\frac Ns}(\R^N)$ such that
	\begin{equation}\label{logFF}
		\bigg|\int_{\R^N}\bigg(\log\frac1{|\cdot|}\ast F(u)\bigg)F(u)\dd x\bigg|<+\infty\,.
	\end{equation}
\end{theorem}
\noindent Next, inspired by some ideas in \cite{BCT}, we rigorously establish the fact that from a solution of \eqref{Chs} one obtains a solution of \eqref{SPs}; in doing this we also investigate regularity and decaying properties of the obtained solutions.
\begin{theorem}[\eqref{Chs}$\,\,\Longrightarrow\,\,$\eqref{SPs}]\label{thm_equiv}
	Assume (V) and ($f_1$)-($f_2$), let $u\in W^{s,\frac Ns}(\R^N)$ be a positive weak solution of the Choquard equation \eqref{Chs} in the sense of Definition \ref{sol_Choquard_log} and define $\phi_u:=I_N\ast F(u)$. Then, $(u,\phi_u)\in W^{s,\frac Ns}(\R^N)\times L_\gamma(\R^N)$ for all $\gamma>0$ is a solution of the Schr\"odinger-Poisson system \eqref{SPs} in the sense of Definition \ref{sol_SP}. Moreover $u\in L^\infty(\R^N)\cap C^{0,\nu}_{\loc}(\R^N)$ for some $\nu\in(0,1)$, and
	\begin{equation}\label{asympt_phi}
		\phi_u(x)=-C_N\|F(u)\|_1\log|x|+o(1)\qquad\mbox{as}\ \ |x|\to+\infty\,.
	\end{equation}
\end{theorem}
\noindent Let us stress the fact that we do not prove that the two problems are equivalent, namely that they have the same set of solutions. Finding a proper functional setting in which this holds is still open, even for the local case \eqref{SP}--\eqref{Ch}, see \cite[Section 2]{BCT}. However, combining Theorems \ref{Thm_Chs} and \ref{thm_equiv} one directly deduces the existence results for \eqref{SPs}, which we next state for the sake of completeness.
\begin{corollary}[Existence for \eqref{SPs}]
	Suppose $V$ is radially symmetric and satisfies (V), and that ($f_1$)-($f_5$) are fulfilled, then \eqref{SPs} possesses a solution $(u,\phi)\in W^{s,\frac Ns}(\R^N)\times L_\gamma(\R^N)$ for all $\gamma>0$ such that:
	\begin{itemize}
		\item[i)] $u\in L^\infty(\R^N)\cap C^{0,\nu}_{\loc}(\R^N)$ for some $\nu\in(0,1)$, is positive, radially symmetric and \eqref{logFF} holds;
		\item[ii)] $\phi=\phi_u:=I_N\ast F(u)$ and the asymptotic behaviour \eqref{asympt_phi} is satisfied.
	\end{itemize}
\end{corollary}

\noindent Finally, in the special case of a constant potential, and assuming some additional regularity of the solutions in order to be able to deal with the pointwise definition of the operator $(-\Delta)_{\frac Ns}^s$, it is possible to show that solutions are radially symmetric, by means of the moving planes technique and exploiting the connection between \eqref{Chs} and \eqref{SPs} established in Theorem \ref{thm_equiv}. Our result extends to system \eqref{SPs} the corresponding results for $\big(s,\frac Ns\big)$-fractional Schr\"odinger equations obtained in \cite{CL}.
\begin{theorem}[Symmetry for \eqref{Chs}]\label{Thm:symmetric}
	Suppose ($f_1$)-($f_2$) are satisfied, and let $u\in C_{loc}^{1,1}(\R^N)\cap W^{s,\frac Ns}(\R^N)$ be a positive solution of \eqref{Chs} with $V(x)=V_0>0$. Then $u$ is radially symmetric around the origin and monotone decreasing.
\end{theorem}

\noindent \textbf{Overview.} After recalling some preliminary results in the next section, in Section \ref{Section_Choq} we prove the existence result for \eqref{Chs}, namely Theorem \ref{Thm_Chs}, by means of a variational approach together with a uniform asymptotic approximation technique. The quite delicate relationship between \eqref{Chs} and \eqref{SPs}, which yields also existence for \eqref{SPs}, is investigated in Section \ref{Section_SP}. The moving planes argument used in the proof of Theorem \ref{Thm:symmetric} follows exactly the one employed in \cite{CLR} for the planar case and thus we omit here the details.

\medskip

\paragraph{\textbf{Notation.}} For $R>0$ and $x_0\in\R^N$ we denote by $B_R(x_0)$ the ball of radius $R$ and center $x_0$. Given a set $\Omega\subset\R^N$, we denote $\Omega^c:=\R^N\setminus\Omega$, and its characteristic function by $\chi_\Omega$. The space of the infinitely differentiable functions which are compactly supported is denoted by $C^\infty_0(\R^N)$. The norm of the Lebesgue spaces $L^p(\R^N)$ with $p\in[1,+\infty]$ is denoted by $\|\cdot\|_p$. The spaces $C^{0,\nu}(\R^N)$ for $\nu\in(0,1)$ are usual spaces of H\"older continuous functions. For $q>0$ we define $q!:=q(q-1)\cdots(q-\lfloor q\rfloor)$, where $\lfloor q\rfloor$ denotes the largest integer strictly less than $q$; if $q>1$ its conjugate H\"older exponent is $q':=\frac q{q-1}$. The symbol $\lesssim$ indicates that an inequality holds up to a multiplicative constant depending only on the structural constants and not on the functions involved. Finally, $o_n(1)$ denotes a vanishing real sequence as $n\to+\infty$. Hereafter, the letter $C$ will be used to denote positive constants which are independent of relevant quantities and whose value may change from line to line.

\section{Preliminaries}\label{Sec_2}
\noindent In this section we summarize some consequences of our assumptions, we recall some basic estimates as well as two well-known results, Lions' compact embedding and the Hardy-Littlewood-Sobolev inequality, which will be used throughout the whole manuscript.

\begin{remark}\label{Rmk_ass} From our set of assumptions can be derived the following: 
	\begin{enumerate}
		\item[(i)] From ($f_1$) and ($f_2$) it is easy to infer that for all $\varepsilon>0$ there exists a constant $C_\varepsilon>0$ such that
		$$f(t)\leq\varepsilon t^{\frac Ns-1}+C_\varepsilon t^{\frac Ns-1}\Phi_{N,s}(\alpha_*|t|^{\frac N{N-s}})\quad\mbox{for all}\ \,t>0\,,$$
		as well as
		$$F(t)\leq\varepsilon t^\frac Ns+C_\varepsilon t^\frac Ns\Phi_{N,s}(\alpha_*|t|^{\frac N{N-s}})\quad\mbox{for all}\ \,t>0\,;$$
		\item[(ii)] Assumption $(f_3)$ implies
		$$F(t)\leq(s-\tau)tf(t)\quad\ \mbox{for any}\ \ t\geq0\,;$$
		\item[(iii)] From ($f_4$) one has that for any $\varepsilon>0$ there exists $M_\varepsilon$ such that for $t\geq M_\varepsilon$
		\begin{equation}\label{Rmk_ass_2_eq}
			F(t)\leq\varepsilon tf(t)\,,
		\end{equation}
		and that there exist $t_0,M_0>0$ such that
		\begin{equation}\label{Rmk_ass_f7}
			F(t)\leq M_0f(t)\quad\ \mbox{for any}\ \ t\geq t_0\,.
		\end{equation}
	\end{enumerate}
	For the proof of (i)-(iii), we refer to \cite[Remark 2.3]{CLR} and \cite[p.2 (1.3)]{CDL}.
\end{remark}

\begin{lemma}\label{Lem_basic_est}
	Let $\mu\in(0,1]$. Then, the following elementary inequality holds 
	$$\frac{t^{-\mu}-1}{\mu}\geq\log\frac1t\qquad\mbox{for all}\ \ t\in(0,1]\,.$$
	Moreover, for all $\nu>\mu$ there exists $C_\nu>0$ such that
	$$
	\frac{t^{-\mu}-1}{\mu}\leq C_\nu t^{-\nu}\qquad\mbox{for all}\ \ t>0\,.
	$$
\end{lemma}

\begin{lemma}[Lemma 2.3, \cite{LY}]\label{estimate_Sani}
	Let $\alpha>0$ and $r>1$. Then for any $\beta>r$ there exists a constant $C_\beta>0$ such that
	\begin{equation*}
		\left(\Phi_{N,s}\left(\alpha|t|^{\frac N{N-s}}\right)\right)^r\leq C_\beta\Phi_{N,s}\left(\alpha\beta|t|^{\frac N{N-s}}\right)\qquad\mbox{for all}\ \,t>0.
	\end{equation*}
\end{lemma}

We denote by $\Wrad^{s,\frac Ns}(\R^N)$ the subspace of radial functions in $W^{s,\frac Ns}(\R^N)$.

\begin{lemma}(see \cite{Lions82})\label{Lem:lions}
	Let $s\in(0,1]$. Then $\Wrad^{s,\frac Ns}(\R^N)$ is compactly embedded in $L^q(\R^N)$ for any $\frac Ns<q<\infty$.
\end{lemma}

\begin{lemma}(Hardy-Littlewood-Sobolev inequality, \cite[Theorem 4.3]{LL})\label{HLS}
	Let $N\geq1$, $\mu\in(0,N)$, and $q,r>1$ with $\tfrac1q+\tfrac\mu N+\tfrac1r=2$. There exists a constant $C=C(N,\mu,q,r)$ such that for all $f\in L^q(\R^N)$ and $h\in L^r(\R^N)$ one has
	$$
	\int_{\R^N}\left(\frac1{|\cdot|^\mu}\ast f\right)\!h\dd x\leq C\|f\|_q\|h\|_r\,.
	$$
\end{lemma}

\section{Existence results for \eqref{Chs} by asymptotic approximation:\\ Proof of Theorem \ref{Thm_Chs}}\label{Section_Choq}
\subsection{The approximating method}
As mentioned in the Introduction, the energy functional formally associated to \eqref{Chs}
\begin{equation*}
	J(u):=\frac sN\|u\|_V^{\frac Ns}-\frac{C_N}2\intN\!\left(\intN\log\frac1{|x-y|}F(u(y))\dd y\right)F(u(x))\dd x
\end{equation*}
is not well-defined on the natural Sobolev-Slobodeckij space $W^{s,\frac Ns}(\R^N)$ because of the presence of the logarithmic convolution term. This prevents the use of standard variational tools. Following \cite{CLR}, we will make use of the approximating method, based on the simple convergence \eqref{key_convergence}, which has been developed in the local case in \cite{LRTZ,CDL}. We set
$$G_\mu(x):=\frac{|x|^{-\mu}-1}\mu,\qquad\mu\in(0,1]\,,\ \ x\in\R^N,$$
and consider the approximating problem
\begin{equation}\label{Choq_approx}
	(-\Delta)^s_{\frac Ns} u+V(x)|u|^{\frac Ns-2}u=C_N\left(G_\mu(\cdot)\ast F(u)\right)f(u)\quad\ \mbox{in}\ \ \R^N,
\end{equation}
whose related energy functional is given by 
\begin{equation*}
	\begin{split}
		J_\mu(u):&=\frac sN\|u\|_V^{\frac Ns}-\frac{C_N}2\intN(G_\mu(\cdot)\ast F(u))F(u)\dd x\,\\
		&=\intN\intN\frac{|u(x)-u(y)|^{\frac Ns}}{|x-y|^{2N}}\dd x\dd y+\intN V(x)|u|^{\frac Ns}\dd x\\
		&\quad+\frac{C_N}{2\mu}\bigg[\intN F(u)\dd x\bigg]^2-\frac{C_N}{2\mu}\intN\intN\frac1{|x-y|^\mu}F(u(x))F(u(y))\dd x\dd y.
	\end{split}
\end{equation*}
\noindent The advantage of this approach is that the power-type singularity in $G_\mu$ can be handled by means of the Hardy-Littlewood-Sobolev inequality (Lemma \ref{HLS}), and it is thus possible to prove under conditions (V) and ($f_0$)-($f_2$) that $J_\mu$ is well-defined and $C^1$ on $W^{s,\frac Ns}(\R^N)$ with Fr\'{e}chet derivative
\begin{equation}\label{J_der}
	\begin{split}
		J'_\mu(u)[\varphi]=&\intN\frac{|u(x)-u(y)|^{\frac Ns-2}(u(x)-u(y))(\varphi(x)-\varphi(y))}{|x-y|^{2N}}+\intN V(x)|u|^{\frac Ns-2}u\varphi\dd x\\
		&-C_N\intN(G_\mu(\cdot)\ast F(u))f(u)\varphi\dd x\,.
	\end{split}
\end{equation}
\noindent We aim at proving first that, for all $\mu$ small, the functional $J_\mu$ possesses a critical point $u_\mu$ of mountain pass type. Being able to estimate the norm in $W^{s,\frac Ns}(\R^N)$ of such critical points uniformly with respect to $\mu$, thanks to a careful estimate on the mountain pass level, we will then pass to the limit as $\mu\to0$ in order to show the existence of a critical point for the original functional $J$, which a posteriori will satisfy \eqref{logFF}.

\noindent Since the potential $V$ is radially symmetric, we will work in the radial subspace $\Wrad^{s,\frac Ns}(\R^N)$, where we can exploit the compactness given by Lemma \ref{Lem:lions}.

\noindent Let us start showing that for all $\mu\in(0,1]$ the approximating functional $J_\mu$ enjoys a mountain pass geometry.
\begin{lemma}\label{Lem:MP1}
	Let $\mu\in(0,1]$ and assume ($f_1$)-($f_3$). Then, there exist constants $\rho,\eta>0$ and $e\in\Wrad^{s,\frac Ns}(\R^N)$ such that:
	\begin{enumerate}
		\item[{\rm(i)}] $\|e\|>\rho$ and $J_\mu(e)<0$;
		\item[{\rm(ii)}] $J_\mu|_{S_\rho}\geq\eta>0$, where $S_\rho=\big\{u\in\Wrad^{s,\frac Ns}(\R^N)\,|\,\|u\|=\rho\big\}$.
	\end{enumerate}
\end{lemma}
\begin{proof}
	{\rm(i)} Take $e_0\in\Wrad^{s,\frac Ns}(\R^N)$ with values in $[0,1]$, $\supp e_0\subset B_\frac14(0)$, and such that $e_0\equiv 1$ in $B_\frac18(0)$. For $t>0$ set
	\begin{equation}\label{eqn:AR-0}
		\Psi(t):=\frac12\left(\intN F(te_0)\dd x\right)^2\!.
	\end{equation}
	By Remark \ref{Rmk_ass} we have
	\begin{equation*}
		\Psi'(t)=\frac1t\left(\intN F(te_0)\right)\left(\intN f(te_0)te_0\right)\geq\frac{2\Psi(t)}{t(s-\tau)}\,,
	\end{equation*}
	hence integrating on $[1,t]$ one gets $\Psi(t)\geq\Psi(1)t^{\tfrac 2{s-\tau}}$. Inserting this in the functional $J_\mu$ and using Lemma \ref{Lem_basic_est} we infer
	\begin{equation*}
		\begin{split}
			J_\mu(te_0)&\leq\frac sNt^{\frac Ns}\|e_0\|_V^{\frac Ns}-\frac{C_N}2\int_{\{|x-y|\leq\tfrac12\}}\log\frac1{|x-y|}F(te_0(y))F(te_0(x))\dd x\dd y\\
			&\leq\frac sNt^{\frac Ns}\|e_0\|_V^{\frac Ns}-\frac{C_N\log 2}2\left(\intN F(te_0)\right)^2\\
			&\leq\frac sNt^{\frac Ns}\|e_0\|_V^{\frac Ns}-\frac{C_N\log 2}2\left(\intN F(e_0)\right)^2t^{\tfrac 2{s-\tau}}.
		\end{split}
	\end{equation*}
	Since $\tau>\left(1-\tfrac2N\right)s$ by ($f_3$), one deduces that $J_\mu(te_0)\to-\infty$ as $t\to+\infty$. Hence one can take $e:=te_0$ for a sufficiently large $t$.
	
	{\rm(ii)} By Lemma \ref{Lem_basic_est} with $\nu=1$ we have
	\begin{equation*}
		J_\mu(u)\geq\frac sN\|u\|_V^{\frac Ns}-\frac{C_1C_N}2\int_{\{|x-y|<1\}}\frac{F(u(x))F(u(y))}{|x-y|}\dd x\dd y\,.
	\end{equation*}
	Applying the Hardy-Littlewood-Sobolev inequality with $\mu=1$ and $q=r=\tfrac{2N}{2N-1}$ we infer
	\begin{equation*}
		\begin{split}
			\int_{\{|x-y|<1\}}\frac{F(u(x))F(u(y))}{|x-y|}\dd x\dd y&\les\bigg[\intN|u|^{\frac{2N^2}{s(2N-1)}}+\left(\intN|u|^\frac{2N^2\theta'}{s(2N-1)}\right)^\frac1{\theta'}\\
			&\quad\cdot\left(\intN\Phi_{N,s}\left(\alpha_*r\|u\|^\frac N{N-s}\left(\frac u{\|u\|}\right)^{\frac N{N-s}}\right)\right)^\frac1\theta\bigg]^\frac{2N-1}N
		\end{split}
	\end{equation*}
	with $r>\tfrac{2N\theta}{2N-1}$ by Lemma \ref{estimate_Sani}. Since $\theta>1$ is arbitrary, it is sufficient to require $\rho<\left(\frac{2N-1}{2N}\right)^\frac{N-s}N$ to be able to use on $S_\rho$ the Moser-Trudinger inequality provided by Theorem \ref{MT_ineqPR} and uniformly bound the last term. To control the first two, instead, we use the continuous embedding given by Lions' Lemma \ref{Lem:lions}. We end up with
	$$J_\mu(u)\ges\|u\|^{\frac Ns}-\|u\|^\frac{2N}s\geq\eta>0$$
	for all $u\in S_\rho$, up to a smaller $\rho$, as desired.
\end{proof}

\noindent The geometry given by Lemma \ref{Lem:MP1} ensures the existence of a PS-sequence at the mountain pass level
$$c_\mu:=\inf_{\gamma\in\Gamma}\max_{t\in[0,1]}J_\mu(\gamma(t)),$$
where
$$\Gamma:=\big\{\gamma\in C\big([0,1],\Wrad^{s,\frac Ns}(\R^N)\big)\,|\,\gamma(0)=0,\gamma(1)=e\big\}$$
for any fixed $\mu\in(0,1)$, see e.g. \cite{Ekeland}. In other words, there exists a sequence $(u_n^\mu)_n\subset\Wrad^{s,\frac Ns}(\R^N)$ such that
\begin{equation}\label{PSsequence}
	J_{\mu}(u_n^\mu)\to c_\mu\quad\text{and}\quad J'_\mu(u_n^\mu)\to0\quad\text{in}\,\,\left(\Wrad^{s,\frac Ns}(\R^N)\right)'
\end{equation}
as $n\to+\infty$. For the sake of a lighter notation, hereafter we set $u_n:=u_n^\mu\,$.

\begin{remark}\label{rem:MP_level}
	Observe that from the proof of Lemma \ref{Lem:MP1} there exist two constants $\underline c,\overline c>0$ independent of $\mu$ such that $\underline c<c_\mu<\overline c$.
\end{remark}

\noindent Next, we show that PS-sequences at level $c_\mu$ are bounded.
\begin{lemma}\label{Lem:c-bounded}
	Assume that (f$_1$)--(f$_4$) hold. Let $(u_n)_n\subset W^{s,\frac Ns}(\R^N)$ be a PS-sequence of $J_\mu$ at level $c_\mu$, then $(u_n)_n$ is bounded in $W^{s,\frac Ns}(\R^N)$ and
	\begin{equation}\label{Lem:estimates_GFF_GFfu}
		\left|\intN[G_\mu(x)\ast F(u_n)]F(u_n)\dd x\right|<C\,,\quad\quad\left|\intN[G_\mu(x)\ast F(u_n)]f(u_n)u_n\dd x\right|<C\,.
	\end{equation}
\end{lemma}
\begin{proof}
	Define $$v_n:=\begin{cases}
		\frac{F(u_n)}{f(u_n)},&  u_n>0\,,\\
		(s-\tau)u_n, & u_n\leq0\,.\\
	\end{cases}$$
	By Remark \ref{Rmk_ass} one has $|v_n|\les|u_n|$ in $\R^N$, hence $v_n\in L^{\frac Ns}(\R^N)$. Following the computations in the proof of \cite[Lemma 5.6]{CLR}, one gets $[v_n]_{s,\frac Ns}\les[u_n]_{s,\frac Ns}$, as well as
	\begin{equation}\label{seminorm_un_vn}
		\intN\intN\frac{\left|u_n(x)-u_n(y)\right|^{\frac Ns-2}\left(u_n(x)-u_n(y)\right)\left(v_n(x)-v_n(y)\right)}{|x-y|^{2N}}\dd x\dd y\leq(s-\tau)[u_n]_{s,\frac Ns}^{\frac Ns}\,.
	\end{equation}
	In particular $v_n\in\Wrad^{s,\frac Ns}(\R^N)$ and we may use it as a test function for $J_\mu'(u_n)$. In light of \eqref{seminorm_un_vn}, combining the two relations in \eqref{PSsequence}, by the fact that in both expressions the term $\intN[G_\mu(x)\ast F(u_n)]F(u_n)\dd x$ appears due to the choice of the test function $v_n$, one obtains
	\begin{equation*}
		(s-\tau)[u_n]_{s,\frac Ns}^{\frac Ns}+(s-\tau)\int_{\{u_n<0\}}V(x)|u_n|^{\frac Ns}+\int_{\{u_n\geq0\}}V(x)|u_n|^{\frac Ns-2}u_n\frac{F(u_n)}{f(u_n)}+2c_\mu-\frac{2N}s\|u\|_V^{\frac Ns}\geq o_n(1)\,.
	\end{equation*}
	Using again Remark \ref{Rmk_ass} to control the third term, the above expression reduces to
	\begin{equation}\label{unif_bound_norm_u_n}
		\left(\tau-\left(1-\frac2N\right)s\right)\|u_n\|_V^{\frac Ns}\leq2c_\mu+o_n(1)\,,
	\end{equation}
	hence a uniform bound in $W^{s,\frac Ns}(\R^N)$. The bounds in \eqref{Lem:estimates_GFF_GFfu} then follow combining \eqref{unif_bound_norm_u_n} and \eqref{PSsequence} with $\varphi=v_n$.
\end{proof}

\begin{remark}\label{Rmk_nonneg_Cerami}
	Arguing as in \cite[Remark 5.7]{CLR}, combining the uniform bound given by Lemma \ref{Lem:c-bounded} and \eqref{PSsequence}, from now on we may assume that PS-sequences at level $c_\mu$ are nonnegative.
\end{remark}

\noindent In order to have a more precise bound on the norm of $(u_n)_n$, in light of \eqref{unif_bound_norm_u_n}, we are lead to a careful estimate of the mountain pass level $c_\mu$. Usually, in presence of nonlinearities with exponential critical growth, this is accomplished combining a growth assumption in the spirit of de Figueiredo-Miyagaki-Ruf \cite{dFMR} with fine estimates by means of a Moser sequence, which detects the critical Moser level in the Trudinger inequality. However, in the fractional setting this way is not feasible, since it is still not known whether the fractional analogue of the Moser sequence proposed in \cite{PR} detects the critical Moser value of the fractional Trudinger inequality, see Theorem \ref{MT_ineqPR}. A possible alternative way, often used in the literature, is to prescribe a strong control of $f$ near zero and with a large constant, the upper bound of which is often not explicit, see e.g. \cite{AF,CMR,AFS,BM1,BM2,ChWa}. We propose instead an approach which combines the choice of a simple fixed test function, which allows explicit computations on its seminorm, and our assumption ($f_5$), in which the growth of $F(\cdot)f(\cdot)$ is bounded from below away from zero by a suitable power multiplied by a constant, on which we have an explicit lower bound. Since it is automatically satisfied by the exponential growth of $f$ at infinity due to ($f_3$)-($f_4$), ($f_5$) is a sort of middle-range assumption, which in most of the cases it is easily verifiable. 

\noindent First, let us introduce our test function and compute its norm. Let $R>0$ and $\bw\in C(\R^N)$ be such that
\begin{equation}\label{bar w}
	\bw(x):=\begin{cases}
		1\quad&\mbox{for}\ |x|\leq\frac R2\,,\\
		2-\frac2R|x|\quad&\mbox{for}\ |x|\in(\frac R2,R)\,,\\
		0\quad&\mbox{for}\ |x|\geq R\,.
	\end{cases}
\end{equation}
\begin{lemma}\label{Lemma_bar w}
	For all $R\in(0,1]$ and $s\in(0,1)$ we have $\bw\in\Wrad^{s,\frac Ns}(\R^N)$ and
	$$\|\bw\|_V^{\frac Ns}\leq\mathfrak J_N(s,R)+\mathfrak K_N(s),$$
	where the constants $\mathfrak J_N(s,R)$ and $\mathfrak K_N(s)$ are given explicitly in \eqref{norm_w} and \eqref{seminorm_w}, respectively.
\end{lemma}
\begin{proof}
	First, by ($V$) one has
	\begin{equation}\label{norm_w_first}
		\intN V(x)|\bw|^{\frac Ns}\dd x\leq\overline{V_0}\left(|B_{\frac R2}(0)|+N\omega_{N-1}\int_{\frac R2}^R\left|2-\frac2R\rho\right|^{\frac Ns}\rho^{N-1}\dd\rho\right).
	\end{equation}
	By the change of variable $y=R-\rho$, we estimate
	\begin{equation*}
		\begin{split}
			\int_{\frac R2}^R&\left|2-\frac2R\rho\right|^{\frac Ns}\rho^{N-1}\dd\rho=\left(\frac2R\right)^{\frac Ns}\int_0^{\frac R2}y^{\frac Ns}(R-y)^{N-1}\dd y\\
			&\leq \left(\frac2R\right)^{\frac Ns}\int_0^{\frac R2}y^{\frac Ns}(R-y)\dd y=\frac{s(N+3s)R^2}{4(N+s)(N+2s)}\,.
		\end{split}
	\end{equation*}
	Combining this with \eqref{norm_w_first}, we obtain
	\begin{equation}\label{norm_w}
		\intN V(x)|\bw|^{\frac Ns}\dd x\leq\overline{V_0}\left(\omega_N\left(\frac R2\right)^N+\omega_{N-1}\frac{Ns(N+3s)R^2}{4(N+s)(N+2s)}\right)=:\mathfrak J_N(s,R)\,.
	\end{equation}
	Since $\bw$ is radially symmetric, according to \cite[Proposition 4.3]{PR}, we have
	\begin{equation}\label{seminorm_w_splitting}
		\begin{split}
			[\bw]_{s,\frac Ns}^\frac Ns&=(N\omega_N)^2\int_0^{+\infty}\int_0^{+\infty}|\bw(r)-\bw(t)|^{\frac Ns}r^{N-1}t^{N-1}\frac{r^2+t^2}{|r^2-t^2|^{N+1}}\dd r\dd t\\
			&=2(N\omega_N)^2\left(\int_0^\frac R2\int_\frac R2^R\ +\ \int_0^\frac R2\int_R^{+\infty}\
			+\ \int_\frac R2^R\int_R^{+\infty}\right)+(N\omega_N)^2\int_\frac R2^R\int_\frac R2^R\\
			&=:2(N\omega_N)^2\left(I_1+I_2+I_3\right)+(N\omega_N)^2I_4
		\end{split}
	\end{equation}
	and let us next estimate $I_i$, $i=1,\dots,4$, separately. Recalling that
	\begin{equation}\label{derivative_PR}
		\frac\dd{\dd r}\left(\frac1N\frac{r^N}{(t^2-r^2)^N}\right)=r^{N-1}\frac{r^2+t^2}{|t^2-r^2|^{N+1}}\,,
	\end{equation}
	we have
	\begin{equation}\label{I_1}
		\begin{split}
			I_1&=\int_{\frac R2}^R\left|\frac2Rt-1\right|^{\frac Ns}t^{N-1}\left(\int_0^{\frac R2}r^{N-1}\frac{r^2+t^2}{|t^2-r^2|^{N+1}}\dd r\right)\dd t\\
			&=\frac1N\left(\frac R2\right)^{N-\frac Ns}\int_{\frac R2}^R\left|t-\frac R2\right|^{\frac Ns}\frac{t^{N-1}}{\left(t-\frac R2\right)^N\left(t+\frac R2\right)^N}\dd t\\
			&\leq\frac{R^{N-1}}N\left(\frac2R\right)^{\frac Ns}\left[\frac{\left(t-\frac R2\right)^{\frac Ns-N+1}}{\frac Ns-N+1}\right]_{t=\tfrac R2}^{t=R}=\frac1{N2^{N-1}\left(\frac Ns-N+1\right)}\,.
		\end{split}
	\end{equation}
	Concerning the second term we get similarly
	\begin{equation}\label{I_2}
		I_2=\frac{R^N}N\int_0^{\frac R2}\frac{r^{N-1}}{(R^2-r^2)^N}\dd r\leq \frac1N\left(\frac R2\right)^{N-1}\int_0^{\frac R2}\frac{\dd r}{(R-r)^N}=\frac{2^{N-1}-1}{2^{N-1}(N-1)N}\,.
	\end{equation}
	The third term can be estimated as
	\begin{equation}\label{I_3}
		\begin{split}
			I_3&=\int_{\frac R2}^R\left|\frac2Rr-1\right|^{\frac Ns}r^{N-1}\left(\int_R^{+\infty}t^{N-1}\frac{r^2+t^2}{(t^2-r^2)^{N+1}}\dd t\right)\dd r\\
			&=\frac1N\left(\frac 2R\right)^{\frac Ns}R^N\int_{\frac R2}^R(r-R)^{\frac Ns-N}\frac{r^{N-1}}{(R+r)^N}\dd r\\
			&\leq\frac{2^{\frac Ns}}NR^{N-1-\frac Ns}\int_0^{\frac R2}y^{\frac Ns-N}\dd y=\frac{2^{N-1}}{N\left(\frac Ns-N+1\right)}\,.
		\end{split}
	\end{equation}
	Finally,
	\begin{equation}\label{I_4}
		\begin{split}
			I_4&=\left(\frac2R\right)^{\frac Ns}\int_{\frac R2}^Rr^{N-1}\left(\int_{\frac R2}^R|r-t|^{\frac Ns}t^{N-1}\frac{r^2+t^2}{|t^2-r^2|^{N+1}}\dd t\right)\dd r\\
			&=\left(\frac2R\right)^{\frac Ns}\int_\frac R2^Rr^{N-1}\bigg[|r-t|^{\frac Ns}\cdot\frac1N\frac{t^N}{|t^2-r^2|^N}\bigg]_{t=\frac R2}^{t=R}\dd r\\
			&\quad+\left(\frac2R\right)^{\frac Ns}\int_\frac R2^R\frac{r^{N-1}}s\int_{\frac R2}^r\bigg[(r-t)^{\frac Ns-1}\cdot\frac{t^N}{(r^2-t^2)^N}\bigg]\dd t\dd r\\
			&\quad-\left(\frac2R\right)^{\frac Ns}\int_\frac R2^R\frac{r^{N-1}}s\int_r^R\bigg[(t-r)^{\frac Ns-1}\cdot\frac{t^N}{(t^2-r^2)^N}\bigg]\dd t\dd r\\
			&=:A_1+A_2-A_3\leq A_1+A_2\,,
		\end{split}
	\end{equation}
	since $A_3\geq0$. Moreover,
	\begin{equation}\label{A_1}
		\begin{split}
			A_1&\leq\left(\frac2R\right)^{\frac Ns}\int_{\frac R2}^R\frac{r^{N-1}}N(R-r)^{\frac Ns-N}\frac{R^N}{(R+r)^N}\dd r\\
			&\leq\left(\frac2R\right)^{\frac Ns}R^{N-1}\int_{\frac R2}^R(R-r)^{\frac Ns-N}\dd r=\frac1{N2^{N-1}\left(\frac Ns-N+1\right)}
		\end{split}
	\end{equation}
	and
	\begin{equation}\label{A_2}
		\begin{split}
			A_2&\leq\left(\frac2R\right)^{\frac Ns}\frac{R^{N-1}}s\int_{\frac R2}^R\int_{\frac R2}^r(r-t)^{\frac Ns-1-N}\dd t\dd r\\
			&\leq\left(\frac2R\right)^{\frac Ns}\frac{R^{N-1}}s\int_{\frac R2}^R\left(r-\frac R2\right)^{\frac Ns-N}\dd r=\frac{2^{N-1}}{s\left(\frac Ns-N+1\right)}\,.
		\end{split}
	\end{equation}
	Eventually, from \eqref{seminorm_w_splitting}-\eqref{A_2} we get
	\begin{equation}\label{seminorm_w}
		[\bw]_{s,\frac Ns}^\frac Ns\leq\frac{N\omega_N^2}{\frac Ns-N+1}\left(\frac 3{2^N}+\frac{(2^N-2)\left(\frac Ns-N+1\right)}{2^N(N-1)}+2^N\left(1+\frac N{2s}\right)\right)=:\mathfrak K_N(s)\,.
	\end{equation}	
\end{proof}

\noindent Thanks to this information on $\bw$ in $\Wrad^{s,\frac Ns}(\R^N)$ we are now able to estimate the mountain pass level via ($f_5$).
\begin{lemma}\label{Lem:MP_level}
	Under ($f_1$)-($f_5$) one has $c_\mu<\frac s{2N}$ for all $\mu\in(0,1]$.
\end{lemma}
\begin{proof}
	Let $\bw\in\Wrad^{s,\frac Ns}(\R^N)$ defined in \eqref{bar w} and $w:=\bw\|\bw\|_V^{-1}$. Then there exists $T>0$ such that $J_\mu(Tw)=\max_{t\geq0}J_\mu(tw)$ and $\tfrac{\dd}{\dd t}|_{t=0}J_\mu(tw)=J_\mu'(Tw)w=0$. Suppose by contradiction that $J_\mu(Tw)\geq\frac s{2N}$, then
	\begin{equation}\label{MPest_1}
		\frac sNT^{\frac Ns}\geq\frac s{2N}+\frac{C_N}2\intN\left(G_\mu(\cdot)\ast F(Tw)\right)F(Tw)\,.
	\end{equation}
	Requiring $R\leq\frac13$, then $G_\mu(x-y)\geq\log\frac1{|x-y|}\geq\log\frac32>0$, namely the convolution term is positive. We therefore obtain from \eqref{MPest_1} that
	\begin{equation}\label{T_below}
		T\geq\left(\frac12\right)^{\frac sN}.
	\end{equation}
	Exploiting the fact that $J_\mu'(Tw)[Tw]=0$, we get
	\begin{equation*}
		\begin{split}
			T^{\frac Ns}&=C_N\intN\left(G_\mu(\cdot)\ast F(Tw)\right)f(Tw)Tw\dd x\\
			&\geq\int_{B_\frac R2(0)}\int_{B_\frac R2(0)}\log\frac1{|x-y|}F(Tw(y))f(Tw(x))Tw(x)\dd x\dd y\\
			&\geq\log3\int_{B_{\frac R2}(0)}F(Tw)\dd y\int_{B_{\frac R2}(0)}f(Tw)Tw\dd x\\
			&\geq\log3\,\bigg(\!\int_{B_{\frac R2}(0)}\sqrt{F(Tw)f(Tw)Tw}\dd x\bigg)^2.
		\end{split}
	\end{equation*}
	Since $\bw\equiv1$ in $B_{\frac R2}(0)$, then
	$$Tw=T\|\bw\|_V^{-1}\geq\left(2^{\frac sN}(\mathfrak J_N(s,R)+\mathfrak K_N(s))\right)^{-\frac sN}$$
	having used \eqref{T_below} and Lemma \ref{Lemma_bar w}. Defining thus
	\begin{equation}\label{T(s)}
		T_N(s,R):=\left(2^{\frac sN}(\mathfrak J_N(s,R)+\mathfrak K_N(s))\right)^{-\frac sN},
	\end{equation}
	we may apply ($f_5$) and get
	\begin{equation*}
		T^{\frac Ns}\geq C_N\log3|B_{\frac R2}(0)|^2\beta(Tw)^{\frac Ns+1},
	\end{equation*}
	which yields, using again \eqref{MPest_1},
	\begin{equation*}
		\beta\leq\beta2^{\frac Ns}T\leq\left(\frac2R\right)^N\frac{2^{\frac Ns}\|\bw\|_V^{\frac Ns+1}}{\omega_N^2C_N\log3}\,.
	\end{equation*}
	Fixing now $R=\tfrac13$ and defining $\mathfrak J_N(s):=\mathfrak J_N(s,\tfrac13)$ and $T_N(s):=T_N(s,\tfrac13)$, from Lemma  \ref{Lemma_bar w} one gets
	\begin{equation}\label{beta0}
		\beta\leq\frac{2^{\frac Ns+N}3^N}{\omega_N^2C_N\log3}\left(\mathfrak J_N(s)+\mathfrak K_N(s)\right)^{\frac sN+1}=:\beta_0,
	\end{equation}
	which finally contradicts the lower bound on $\beta$ in ($f_5$). This proves that $J_\mu(Tw)<\frac s{2N}$ and the conclusion follows from the definition of $c_\mu$.
\end{proof}

\noindent Combining \eqref{unif_bound_norm_u_n} and Lemma \ref{Lem:MP_level} we immediately obtain an estimate on the norm of $(u_n)_n$ \textit{which is uniform with respect to} $\mu$:
\begin{equation}\label{unif_bound_norm_u_n_UNIFORMmu}
	\|u_n\|_V^{\frac Ns}\leq\frac{\frac sN+o_n(1)}{\tau-\left(1-\frac2N\right)s}<\frac s{\tau-\left(1-\frac2N\right)s}\,.
\end{equation}

\noindent In order to get a critical point for the functional $J_\mu$, we need first an integrability result for $F(u_n)$ in Lebesgue spaces.
\begin{lemma}\label{Lem:integral-F}
	Suppose ($f_1$)--($f_5$) hold and let $(u_n)_n$ be a PS-sequence of $J_\mu$ at level $c_\mu$. Then there exists $C>0$ independent of $n$ and $\mu$ such that
	\begin{equation}\label{improved_int_F}
		\intN f(u_n)u_n\dd x\leq C \quad \text { and } \quad \intN F(u_n)^\kappa\dd x\leq C
	\end{equation}
	for any $\kappa\in\big[1,\tfrac1{\gamma_N(s,\tau)}\big)$, where $\gamma_N(s,\tau)\in(0,1)$ is a constant depending just on $N$, $s$, and $\tau$.
\end{lemma}
\begin{proof}
	Let us introduce the following auxiliary function
	$$H_N(t):=t-\frac N{2s}\frac{F(t)}{f(t)}\quad\mbox{for}\ \,t \geq 0\,,$$
	and define $v_n:=H_N\left(u_n\right)$. Similarly to the proof of Lemma \ref{Lem:c-bounded}, one has $v_n\in W^{s,\frac Ns}(\R^N)$. We first \textit{claim} that there exists $\gamma_N(s,\tau)\in(0,1)$ depending just on $N$, $s$, and $\tau$, such that
	\begin{equation}\label{gamma_Nstau}
		\left\|v_n\right\|_V^{\frac Ns}\leq\gamma_N(s,\tau)<1
	\end{equation}
	for $n$ large enough. If so, then by $(f_4)$ for any $\varepsilon>0$ there exists $t_\varepsilon>0$ such that
	$$\frac{F(t)}{f(t)}\leq\varepsilon(t-t_\varepsilon)+\frac{F\left(t_\varepsilon\right)}{f\left(t_\varepsilon\right)}\quad\ \text{for all}\ \ t\geq t_\varepsilon\,.$$
	Hence,
	\begin{equation*}
		v_n=H_N\left(u_n\right)=\left(1-\frac{\varepsilon N}{2s}\right)\left(u_n-t_\varepsilon\right)+t_\varepsilon-\frac N{2s}\frac{F\left(t_\varepsilon\right)}{f\left(t_\varepsilon\right)}\geq\left(1-\frac{\varepsilon N}{2s}\right)\left(u_n-t_\varepsilon\right),
	\end{equation*}
	having applied $(f_3)$ in the last inequality. This implies
	\begin{equation}\label{eqn:integ4+}
		u_n(x)\leq t_\varepsilon+\frac{v_n(x)}{1-\bar\varepsilon}\quad\mbox{for all}\ \,\,x\in\R^N,
	\end{equation}
	where $\bar\varepsilon:=\frac{\varepsilon N}{2s}$.
	Hence, by $\left(f_1\right)$-$\left(f_2\right)$ we deduce that for any given $\varepsilon>0$, there exists $C_\varepsilon>0$ such that
	\begin{equation}\label{eqn:integ5}
		\begin{aligned}
			\intN F(u_n)^\kappa\dd x&\leq\int_{\{|u_n|<t_\varepsilon\}}F(u_n)^\kappa\dd x+\int_{\{|u_n|\geq t_\varepsilon\}}F(u_n)^\kappa\dd x \\
			&\leq C_\varepsilon\left\|u_n\right\|^{\frac Ns\kappa}+C_\varepsilon\int_{\{|u_n|\geq t_\varepsilon\}}\Phi_{N, s}\left(\alpha_*\kappa\theta\left(t_\varepsilon+\frac{v_n}{1-\bar\varepsilon}\right)^{\frac N{N-s}}\right)\dd x\\
			&\leq C_\varepsilon\left\|u_n\right\|^{\frac Ns\kappa}+C_\varepsilon\intN\Phi_{N, s}\left(\alpha_*\kappa\theta(1+\bar\varepsilon)\left(\frac{v_n}{1-\bar\varepsilon}\right)^{\frac N{N-s}}\right)\dd x\\
		\end{aligned}
	\end{equation}
	with $\theta>1$ by Lemma \ref{estimate_Sani}. By \eqref{gamma_Nstau} there exists $\sigma>0$ such that $\left\|v_n\right\|_V^\frac Ns\leq\gamma_N(s,\tau)+\sigma<1$ for $n$ large enough, therefore
	\begin{equation}\label{End_proof_Lemma_kappa}
		\frac{\kappa\theta(1+\bar\varepsilon)}{(1-\bar\varepsilon)^{\frac N{N-s}}}\left\|v_n\right\|_V^{\frac N{N-s}}\leq\frac{\kappa\theta(1+\bar\varepsilon)}{(1-\bar\varepsilon)^{\frac N{N-s}}}\left(\gamma_N(s,\tau)+\sigma\right)<1
	\end{equation}
	for $\kappa\in\big[1,\tfrac1{\gamma_N(s,\tau)}\big)$, for suitable choices of for $\varepsilon>0$ and $\theta>1$ small enough. Thanks to Theorem \ref{MT_ineqPR} with $\lambda=\min\{1,\underline{V_0}^{\frac sN}\}$, this yields the second inequality in \eqref{improved_int_F}. Similarly, but more easily, one can also prove the first inequality.
	
\noindent 	The remaining part of the proof is devoted to justifying of the claim \eqref{gamma_Nstau} above. Combining the two relations in \eqref{PSsequence} as in Lemma \ref{Lem:estimates_GFF_GFfu}, one has
	\begin{equation}\label{eqn:integ3-bis}
		\begin{split}
			&\intN\intN\frac N{2s}\frac{\left|u_n(x)-u_n(y)\right|^{\frac Ns-2}\left(u_n(x)-u_n(y)\right)}{|x-y|^{2N}}\left(\frac{F\left(u_n(x)\right)}{f\left(u_n(x)\right)}-\frac{F\left(u_n(y)\right)}{f\left(u_n(y)\right)}\right)\dd x\dd y\\
			&\quad+\frac N{2s}\intN\left|u_n\right|^{\frac Ns-2}u_n\frac{F\left(u_n\right)}{f\left(u_n\right)}\dd x-\|u_n\|_V^{\frac Ns}-\frac Nsc_\mu=o_n(1)\,.
		\end{split}
	\end{equation}
	Therefore,
	\begin{equation}\label{eqn:integ3}
		\begin{split}
		\|v_n\|_V^{\frac Ns}&=\intN\intN\frac{\left|H_N\left(u_n(x)\right)-H_N\left(u_n(y)\right)\right|^{\frac Ns}}{|x-y|^{2N}}\dd x\dd y+\intN V(x)\left|H_N\left(u_n\right)\right|^{\frac Ns}\!\dd x\\
		&=\intN\intN \frac{\left|H_N\!\left(u_n(x)\right)-H_N\!\left(u_n(y)\right)\right|^{\frac Ns}}{|x-y|^{2N}}\dd x\dd y\\
		&\quad+\intN\!\intN\!\frac{\frac N{2s}\left|u_n(x)-u_n(y)\right|^{\frac Ns-2}\!\left(u_n(x)-u_n(y)\right)\!\left(\frac{F\left(u_n(x)\right)}{f\left(u_n(x)\right)}-\frac{F\left(u_n(y)\right)}{f\left(u_n(y)\right)}\right)}{|x-y|^{2N}}\dd x\dd y\\
		&\quad+\intN V(x)\left(|H_N(u_n)|^{\frac Ns}\!+\frac N{2s}\intN|u_n|^{\frac Ns-2}u_n\frac{F(u_n)}{f(u_n)}\right)\!\!\dd x\\
		&\quad+\intN\frac Nsc_\mu-[u_n]_{s,\frac Ns}^{\frac Ns}-V(x)|u_n|^{\frac Ns}\dd x+o_n(1)\,.
	\end{split}
	\end{equation}
	Let us denote by $Z_N(u_n)$ the first term in the right-hand side. Using the mean value theorem we infer
	\begin{equation*}
		Z_N(u_n)=\intN\!\intN\!\left[\left|1-\tfrac N{2s}\left(1-\tfrac{Ff'}{f^2}\left(\theta(x,y)\right)\right)\right|^{\tfrac Ns}\!\!+\tfrac N{2s}\left(1-\tfrac{Ff'}{f^2}\left(\theta(x,y)\right)\right)\right]\!\frac{|u_n(x)-u_n(y)|^{\tfrac Ns}}{|x-y|^{2N}}\dd x\dd y
	\end{equation*}
	for some $\theta(x,y)\in\left(\min\left\{u_n(x),u_n(y)\right\},\max\left\{u_n(x),u_n(y)\right\}\right)$. Splitting $\R^{2N}=A\cup(\R^{2N}\setminus A)$ and $Z_N(u_n)=Z_A(u_n)+Z_{\R^{2N}\setminus A}(u_n)$ accordingly, we estimate the two terms separately. Since
	$$\left|1-\tfrac N{2s}\left(1-\tfrac{Ff'}{f^2}(\theta)\right)\right|^{\frac Ns}\leq1-\tfrac N{2s}\left(1-\tfrac{Ff'}{f^2}(\theta)\right),$$
	one infers
	\begin{equation}\label{I_R4-A}
		Z_{\R^{2N}\setminus A}(u_n)\leq\int_{\R^{2N}\setminus A}\frac{\left|u_n(x)-u_n(y)\right|^{\frac Ns}}{|x-y|^{2N}}\dd x\dd y\,.
	\end{equation}
	On the other hand, applying \cite[Lemma 5.9]{CLR} with $q=\tfrac Ns$ and $W=\frac N{2s}\left(\tfrac{Ff'}{f^2}(\theta)-1\right)>0$,
	\begin{equation*}
		\begin{split}
			Z_A(u_n)-\int_A\frac{\left|u_n(x)-u_n(y)\right|^{\frac Ns}}{|x-y|^{2N}}\dd x\dd y&\leq\frac{\tfrac Ns!}{\tfrac Ns-\lfloor\tfrac Ns\rfloor}\int_A\frac W{1-W}\frac{|u_n(x)-u_n(y)|^{\frac Ns}}{|x-y|^{2N}}\dd x\dd y\,.
		\end{split}
	\end{equation*}
	Imposing
	\begin{equation}\label{def_mu_1}
		\frac{F(t)f'(t)}{f^2(t)}<1+\frac sN
	\end{equation}
	one finds $\|W\|_\infty<\tfrac12$, from which
	\begin{equation*}
		\begin{split}
			Z_A(u_n)-\int_A\frac{\left|u_n(x)-u_n(y)\right|^{\frac Ns}}{|x-y|^{2N}}&\dd x\dd y\leq\frac{\tfrac Ns!}{\tfrac Ns-\lfloor\tfrac Ns\rfloor}\int_A\frac W{1-W}\frac{|u_n(x)-u_n(y)|^{\frac Ns}}{|x-y|^{2N}}\dd x\dd y\\
			&\leq\frac{2\tfrac Ns!}{\tfrac Ns-\lfloor\tfrac Ns\rfloor}\|W\|_\infty\intN\frac{|u_n(x)-u_n(y)|^{\frac Ns}}{|x-y|^{2N}}\dd x\dd y\\
			&\leq\frac{2s\tfrac Ns!}{\left(\tfrac Ns-\lfloor\tfrac Ns\rfloor\right)\left(\tau-\left(1-\frac2N\right)s\right)}\|W\|_\infty\,,
		\end{split}
	\end{equation*}
	having used \eqref{unif_bound_norm_u_n_UNIFORMmu}. This, combined with \eqref{I_R4-A}, yields
	\begin{equation}\label{seminorm_HN}
		Z(u_n)-[u_n]_{s,\frac Ns}^{\frac Ns}\leq\frac{2s\tfrac Ns!\|W\|_\infty}{\left(\tfrac Ns-\lfloor\tfrac Ns\rfloor\right)\left(\tau-\left(1-\frac2N\right)s\right)}\,.
	\end{equation}
	Moreover, since $s-\tau<\frac2Ns$, from ($f_3$) one has
	\begin{equation}\label{norm_HN}
		\begin{split}
			\intN V(x)|H_N(u_n)|^{\frac Ns}\dd x&=\intN V(x)|u_n|^{\frac Ns}\left|1-\frac N{2s}\frac{F(u_n)}{f(u_n)u_n}\right|^{\frac Ns}\!\dd x\\
			&\leq-\frac N{2s}\intN V(x)|u_n|^{\frac Ns}\frac{F(u_n)}{f(u_n)u_n}\dd x\,.
		\end{split}
	\end{equation}
	Combining \eqref{norm_HN} and \eqref{seminorm_HN} with \eqref{eqn:integ3}, one finally gets
	\begin{equation*}
		\|v_n\|_V^{\frac Ns}\leq\frac{2s\tfrac Ns!\|W\|_\infty}{\left(\tfrac Ns-\lfloor\tfrac Ns\rfloor\right)\left(\tau-\left(1-\frac2N\right)s\right)}+\frac Nsc_\mu+o_n(1)\,.
	\end{equation*}
	Since $c_\mu<\frac s{2N}$ by Lemma \ref{Lem:MP_level}, and recalling the definition of $W$, we obtain \eqref{gamma_Nstau}, provided
	\begin{equation}\label{def_mu_2}
		\left\|\frac{Ff'}{f^2}\right\|_\infty<1+\frac{\tfrac Ns-\lfloor\tfrac Ns\rfloor}{\frac Ns\frac Ns!}\left(\frac{\tau\left(1-\frac2N\right)s}{2s}\right).
	\end{equation}
	Defining now
	\begin{equation}\label{def_mu}
		\mu_N(s,\tau):=\min\left\{\frac sN,\,\frac{\tfrac Ns-\lfloor\tfrac Ns\rfloor}{\frac Ns\frac Ns!}\left(\frac{\tau\left(1-\frac2N\right)s}{2s}\right)\right\},
	\end{equation}
	assumption ($f_3$) guarantees that both conditions \eqref{def_mu_1} and \eqref{def_mu_2} are satisfied and the proof is concluded.
\end{proof}

\noindent We are finally in a position to show the existence of a critical point of mountain pass type for the approximating functional $J_\mu$.
\begin{proposition}\label{Pro:u_mu}
	Assume that ($f_1$)--($f_5$) hold. For all $\mu\in(0,1)$ sufficiently small there exists a positive $u_\mu\in\Wrad^{s,\frac Ns}(\R^N)$ such that $J_\mu'(u_\mu)=0$.
\end{proposition}
\begin{proof}
\noindent 	Having in hands the results obtained so far in this section, the proof follows the lines of \cite[Proposition 5.10]{CLR}. Let us just resume here the main steps.
	
\noindent 	By Lemma \ref{Lem:MP1} and Remark \ref{Rmk_nonneg_Cerami}, for all $\mu\in(0,1)$ there exists a nonnegative PS-sequence $(u_n^\mu)_n\subset\Wrad^{s,\frac Ns}(\R^N)$ of $J_\mu$ at level $c_\mu$, which is bounded in $W^{s,\frac Ns}(\R^N)$ by Lemma \ref{Lem:c-bounded}. Hence there exists a nonnegative $u_\mu\in\Wrad^{s,\frac Ns}(\R^N)$ such that, up to a subsequence,
	\begin{equation}\label{eqn:fun3-0}
		\aligned
		u_n^\mu\rightharpoonup u_\mu\ \quad&\text{in}\,\,W^{s,\frac Ns}(\R^N),\\
		u_n^\mu\to u_\mu\ \quad&\text{a.e. in}\,\, \R^N\ \,\mbox{and in}\,\, L^p(\R^N)\ \,\mbox{for all}\ \,p\in\big(\tfrac Ns,+\infty\big).
		\endaligned
	\end{equation}
	First, we prove that
	\begin{equation}\label{conv_F}
		\intN F(u_n^\mu)\dd x\to\intN F(u_\mu)\dd x\,.
	\end{equation}
	Indeed, by the mean value theorem, there exists $\tau_n(x)\in(0,1)$ such that
	\begin{equation*}
		\begin{split}
			\intN&|F(u_n^\mu)-F(u_\mu)|\dd x
			=\intN|f\left(u_\mu+\tau_n(x)(u_n^\mu-u_\mu)\right)(u_n^\mu-u_\mu)|\dd x\\
			&\leq\varepsilon\intN(|u_n^\mu|+|u_\mu|)^{\frac Ns-1}|u_n^\mu-u_\mu|\dd x\\
			&\quad+C_\varepsilon\intN(|u_n^\mu|+|u_\mu|)^{\frac Ns-1}\Phi_{N,s}\left(\alpha_*|u_\mu+\tau_n(x)(u_n^\mu-u_\mu)|^{\frac Ns}\right)|u_n^\mu-u_\mu|\dd x\\
			&\leq\varepsilon\left(\|u_n^\mu\|_\frac Ns+\|u_\mu\|_\frac Ns\right)\|u_n^\mu-u_\mu\|_\frac Ns+C_\varepsilon\int_{\{u_n^\mu>u_\mu\}}|u_n^\mu|^{\frac Ns-1}\Phi_{N,s}(\alpha_*|u_n^\mu|^{\frac Ns})|u_n^\mu-u_\mu|\\
			&\quad+C_\varepsilon\int_{\{u_n^\mu\leq u_\mu\}}|u_\mu|^{\frac Ns-1}\Phi_{N,s}(\alpha_*|u_n|^{\frac Ns})|u_n^\mu-u_\mu|\\
			&=:C\varepsilon+S_1+S_2
		\end{split}
	\end{equation*}
	by \eqref{unif_bound_norm_u_n_UNIFORMmu}. We have
	\begin{equation*}
		S_1\les\left(\intN|u_n^\mu|^{(\frac Ns-1)\theta'\eta}\right)^\frac1{\theta'\eta}\!\left(\intN\Phi_{N,s}\left(\alpha_*r\left|t_\varepsilon+\frac{v_n}{1-\bar\varepsilon}\right|^{\frac Ns}\right)\right)^\frac1{\theta'\eta'}\!\left(\intN |u_n^\mu-u_\mu|^\theta\right)^\frac1\theta=o_n(1),
	\end{equation*}
	with $r>\theta'\eta'$, by choosing $\theta>\frac Ns$ and $\left(\tfrac Ns-1\right)\theta'\eta>\frac Ns$. The uniform boundedness of the second term can be proved arguing as in the proof of Lemma \ref{Lem:integral-F}. The term $S_2$ can be similarly handled.
	
	\noindent Next, we show by means of Lemma \ref{Lem:integral-F} that there exists $C>0$ independent of $n$ and $\mu\in(0,\mu_0)$ such that for all $x\in\R^N$ one has
	\begin{equation}\label{CLAIM}
		\intN\frac{F(u_n(y))}{|x-y|^\mu}\dd y\leq C\,.
	\end{equation}
	Indeed,
	\begin{equation*}
		\begin{split}
			\intN\frac{F(u_n(y))}{|x-y|^\mu}\dd y&=\int_{\{|x-y|\geq1\}}\frac{F(u_n^\mu(y))}{|x-y|^\mu}\dd y+\int_{\{|x-y|<1\}}\frac{F(u_n^\mu(y))}{|x-y|^\mu}\dd y\\
			&\leq\intN F(u_n^\mu)+\left(\int_{\{|x-y|<1\}}\frac{\dd y}{|x-y|^{\mu q}}\right)^\frac1q\left(\intN F(u_n^\mu)^{q'}\right)^\frac1{q'}.
		\end{split}
	\end{equation*}
	Choosing $q$ sufficiently large so that $q'\in\left(1,\tfrac1{\gamma_N(s,\tau)}\right)$, the first and last terms are bounded by Lemma \ref{Lem:integral-F}; then, it is sufficient to choose $\mu_0$ sufficiently small so that $|\cdot|^{\mu q}\in L^1(B_1(0))$ for all $\mu\in(0,\mu_0)$. As a consequence, one deduces that
	\begin{equation}\label{F_conv_Riesz_mu}
		\intN\left(\frac1{|\cdot|^\mu}\ast F(u_n^\mu)\right)F(u_n^\mu)\dd x\to\intN\left(\frac1{|\cdot|^\mu}\ast F(u_\mu)\right)F(u_\mu)\dd x\,.
	\end{equation}
	Indeed, first by $J_\mu'(u_n)[u_n]=o_n(1)$ and Lemmas \ref{Lem:c-bounded} and \ref{Lem:integral-F} one infers that
	\begin{equation*}
		\intN\left(\frac1{|\cdot|^\mu}\ast F(u_n)\right)f(u_n)u_n\leq C
	\end{equation*}
	uniformly with respect to both $n$ and $\mu\in(0,1)$. This implies that the queues of both integrals in \eqref{F_conv_Riesz_mu} are small, by using \eqref{Rmk_ass_f7}. On the other hand, the convergence in the interior part of the integrals in \eqref{F_conv_Riesz_mu} can be proved via the dominated convergence theorem relying on \eqref{CLAIM} and \cite[Lemma 2.1]{dFMR}; see \cite[Lemma 3.6]{R} for more details.
	
	\noindent In a similar way one can also prove that for all $\varphi\in\Wrad^{s,\frac Ns}(\R^N)$ one has
	\begin{equation}\label{conv_f}
		\intN f(u_n^\mu)\varphi\dd x\to\intN f(u_\mu)\varphi\dd x
	\end{equation}
	and
	\begin{equation}\label{f_conv_Riesz_mu}
		\intN\left(\frac1{|\cdot|^\mu}\ast F(u_n^\mu)\right)f(u_n^\mu)\varphi\dd x\to\intN\left(\frac1{|\cdot|^\mu}\ast F(u_\mu)\right)f(u_\mu)\varphi\dd x\,.
	\end{equation}
	\noindent Combining then \eqref{conv_F}, \eqref{conv_f}, \eqref{f_conv_Riesz_mu}, and the weak convergence $u_n^\mu\rightharpoonup u_\mu$ in $W^{s,\frac Ns}(\R^N)$, we deduce that $u_\mu$ is a critical point of $J_\mu$. By exploiting the monotonicity of the operator $(-\Delta)_{\frac Ns}^s$, one can also get $u_n^\mu\to u_\mu$ strongly in $W^{s,\frac Ns}(\R^N)$, which implies, by the continuity of $J_\mu$, that $J_\mu(u_n^\mu)\to J_\mu(u_\mu)$ and so $J_\mu(u_\mu)=c_\mu\geq\underline c$ by Remark \ref{rem:MP_level}. This readily implies that $u_\mu\neq0$. Its positivity follows by the strong maximum principle for the $p-$fractional Laplacian \cite[Theorem 1.4]{Del17}.
\end{proof}

\noindent By Proposition \ref{Pro:u_mu} we obtained a positive critical point for $J_\mu$ for all $\mu$ small, say $\mu\in(0,\bar\mu)$. Note that all relevant estimates obtained so far, in particular \eqref{unif_bound_norm_u_n_UNIFORMmu}, but also \eqref{improved_int_F}, \eqref{CLAIM}, as well as the one for the mountain pass level of Lemma \ref{Lem:MP_level}, are independent of $\mu$. Therefore, thanks to the convergences \eqref{eqn:fun3-0}, \eqref{conv_F}, and $u_n^\mu\to u_\mu$ in $W^{s,\frac Ns}(\R^N)$, one can further get the existence of $u_0\in\Wrad^{s,\frac Ns}(\R^N)$ such that
\begin{equation}\label{conv_umu_u0}
	\aligned
	u_\mu\rightharpoonup u_0\ \quad&\text{in}\,\,W^{s,\frac Ns}(\R^N),\\
	u_\mu\to u_0\ \quad&\text{a.e. in}\,\, \R^N\ \,\mbox{and in}\,\, L^p(\R^N)\ \,\mbox{for all}\ \,p\in\big(\tfrac Ns,+\infty\big).
	\endaligned
\end{equation}
In light of Lemma \ref{Lem:integral-F}, since the constants are uniform with respect to $\mu$, it is easy to infer that
\begin{equation}\label{improved_integrability_F_mu}
	\intN F(u_\mu)\dd x\leq C\quad\mbox{and}\quad\int_{\{|x-y|\leq 1\}}\frac{F(u_\mu(y))}{|x-y|^\frac{4(\omega-1)}{3\omega}}\dd y\leq C
\end{equation}
for $\omega\in\big[1,\tfrac1{\gamma_N(s,\tau)}\big)$. Here, to prove the second inequality from the first, one argues as for \eqref{CLAIM}. Actually, closely following the strategy in \cite[Lemma 5.11]{CLR}, which is based on \eqref{eqn:integ4+}-\eqref{eqn:integ5}, one may also prove that
\begin{equation}\label{Lem:eq_infty}
	\int_{\{|x-y|\leq 1\}}\frac{F(u_\mu(y))}{|x-y|^\frac{4(\omega-1)}{3\omega}}\dd y\to 0\quad\ \ \mbox{as}\ \ |x|\to+\infty\,,
\end{equation}
uniformly for $\mu\in(0,\bar\mu)$. Before we proceed, we need some regularity properties for the solution sequence $(u_\mu)_\mu$.

\begin{lemma}\label{Lem:Regularity_u_mu}
	Let $\mu\in(0,\bar\mu)$ and $u_\mu\in W^{s,\frac Ns}(\R^N)$ be a weak solution of \eqref{Choq_approx}, then there exists $C>0$ independent of $\mu$ such that $\|u_\mu\|_\infty\leq C$. Furthermore, there exists $R>0$ such that for $|x|\geq R$,
	\begin{equation}\label{decay_estimate}
		|u_\mu(x)|\les\frac1{1+|x|^\frac{s(2N+3)}{2(N-s)}}
	\end{equation}
	uniformly for $\mu\in(0,\frac{4(\omega-1)}{3\omega})$.
\end{lemma}
\begin{proof}
	One can prove the boundedness in $L^\infty(\R^N)$ of the sequence $(u_\mu)_\mu$ following the approach of \cite{Bisci2022} for the $\big(s,\frac Ns\big)$-fractional Schr\"odinger equation, which exploits a Nash-Moser iteration technique. It is easy to adapt to our setting the computations in \cite[Lemma 5.12]{CLR} which are detailed for the planar case $N=2$.
	
	\noindent In order to prove the decay estimate, we start by showing that there exist $C,C_0>0$ such that for all $x\in\R^N$ one has
	\begin{equation}\label{gmu_estimate}
		\intN G_\mu(x-y)F(u_\mu(y))\dd y\leq\frac C\mu\left(\left(\frac{|x|}2\right)^{-\mu}-1\right)+C_0\,.
	\end{equation}
	Indeed, first
	\begin{equation}\label{gmu_estimate_ysmall}
		\int_{\left\{|y|\leq\frac{|x|}2\right\}}G_\mu(x-y)F(u_\mu(y))\dd y\leq\frac1\mu\left(\left(\tfrac{|x|}2\right)^{-\mu}-1\right)\intN F(u_\mu)\dd y\leq\frac C\mu\left(\left(\tfrac{|x|}2\right)^{-\mu}-1\right)
	\end{equation}
	by \eqref{improved_integrability_F_mu}. On the other hand,
	\begin{equation}\label{gmu_estimate_ybig1}
		\int_{\left\{|y|\geq\frac{|x|}2\right\}\cap B_1(x)^c}G_\mu(x-y)F(u_\mu(y))\dd y\leq0\,,
	\end{equation}
	since in $B_1(0)^c$ one has $G_\mu(z)\leq0$, while
	\begin{equation}\label{gmu_estimate_ybig2}
		\int_{\left\{|y|\geq\frac{|x|}2\right\}\cap B_1(x)}G_\mu(x-y)F(u_\mu(y))\dd y\les\int_{B_1(0)}G_\mu(z)\dd z=C_0<+\infty\,.
	\end{equation}
	The estimate \eqref{gmu_estimate} follows by combining \eqref{gmu_estimate_ysmall}-\eqref{gmu_estimate_ybig2}. Since the right-hand side of \eqref{gmu_estimate} pointwisely converges to $-C\log\frac{|x|}2+C_0$, there exists $R_1>0$ independent of $\mu$ small such that for all $\varphi\in W^{s,\frac Ns}(\R^N)$ such that $\varphi\geq0$ and $\supp\varphi\subset B_{R_1}(0)^c$ one has
	\begin{equation}\label{bilinear_form_negative}
		\begin{split}
			\intN\intN&\frac{|u_\mu(x)-u_\mu(y)|^{\frac Ns-2}(u_\mu(x)-u_\mu(y))(\varphi(x)-\varphi(y))}{|x-y|^{2N}}+\intN V(x)|u_\mu|^{\frac Ns-2}u_\mu\varphi\dd x\\
			&=C_N\intN\left(\intN G_\mu(x-y)F(u_\mu(y))\dd y\right)f(u_\mu(x))\varphi(x)\dd x<0
		\end{split}
	\end{equation}
	as $|x|>R_1$. Defining now
	\begin{equation}\label{w}
		w(x)=\frac1{1+|x|^a}\quad\mbox{with}\ \,a=\frac{s(2N+3)}{2(N-s)}\,,
	\end{equation}
	since $a>\tfrac{Ns}{N-s}$, by \cite[Lemma 7.1]{Del20} there exists $R_2>0$ such that
	\begin{equation}\label{w_fractional}
		(-\Delta)_{\frac Ns}^sw(x)\les\frac1{|x|^{2N}}\quad\mbox{for all}\ \,|x|>R_2\,,
	\end{equation}
	therefore
	\begin{equation}\label{w_operator}
		(-\Delta)_{\frac Ns}^sw(x)+V(x)|w|^{\frac Ns-2}w\geq0\quad\mbox{for all}\ \,|x|>\widetilde R_2
	\end{equation}
	for a suitable $\widetilde R_2\geq R_2$. Since $\|u_\mu\|_\infty\leq C$ for all $\mu\in(0,\bar\mu)$, there exists $C_1>0$ such that
	\begin{equation*}
		\psi(x):=u_\mu(x)-C_1w(x)\leq0\quad\quad\text{for}\quad |x|=R_3:=\max\{R_1,\widetilde R_2\}.
	\end{equation*}
	Define now $\widetilde\psi:=\psi^+\chi_{\R^N\setminus B_{R_3}(0)}$, which belongs to $\Wrad^{s,\frac Ns}(\R^N)$ and has $\supp\widetilde\psi\subset B_{R_1}(0)^c$ by construction, use it as a test function in \eqref{bilinear_form_negative} and, arguing as in the final part of the proof of \cite[Lemma 3.4]{CLR}, one finds $\widetilde\psi>0$ in $\R^N\setminus B_{R_3}(0)$, that is, the decay estimate in \eqref{decay_estimate}.
\end{proof}

\noindent Now we are in a position to prove the existence of a nontrivial critical point for the original functional $J$. In view of the results obtained so far, the strategy closely follows the one used in \cite{CLR}, see also \cite{R}. Here we report only the main steps.

\begin{proof}[Proof of Theorem \ref{Thm_Chs}]
	Since $J_\mu'(u_\mu)=0$, \eqref{bilinear_form_negative} holds for any test function $\varphi\in\Wrad^{s,\frac Ns}(\R^N)$. The left hand-side converges to the respective for $u_0$ thanks to the weak convergence in \eqref{conv_umu_u0}. The right-hand side will be handled via dominated convergence theorem, by using Lemma \ref{Lem:Regularity_u_mu}. Let us split the argument according to $|x-y|\gtreqqless1$.
	
\noindent 	If $|x-y|\leq1$, we define
	\begin{equation*}
		\begin{split}
			\underline{L_\mu}(x,y)&=G_\mu(x-y)F(u_\mu(y))f(u_\mu(x))\varphi(x)\chi_{\{|x-y|\leq1\}}(x,y)\\
			&\leq C_\omega|x-y|^{-\frac{4(\omega-1)}{3\omega}}F(u_\mu(y))f(u_\mu(x))\varphi(x)\chi_{\{|x-y|\leq1\}}(x,y)=:h(u_\mu)(x,y)
		\end{split}
	\end{equation*}
	and, thanks to \eqref{improved_integrability_F_mu} and the fact that $\intN f(u_\mu)\varphi\dd x\leq C$ uniform with respect to $\mu$, one finds that $h(u_\mu)$ is uniformly bounded in $L^1(\R^N)$. Since in addition $u_\mu\to u_0$ a.e. in $\R^N$, by \cite[Lemma 2.1]{dFMR} one gets $h(u_\mu)\to h(u_0)$, and in turns
	\begin{equation}\label{xysmall}
		\underline{L_\mu}(x,y)\to\log\frac1{|x-y|}F(u_0(y))f(u_0(x))\varphi(x)\chi_{\{|x-y|\leq1\}}(x,y)\quad\mbox{in}\ \, L^1(\R^N)\,.
	\end{equation}
	
	\noindent On the other hand, if $|x-y|>1$ one has $0\geq G_\mu(x-y)=-|x-y|^{\tau\mu}\log|x-y|$ for some $\tau=\tau(x,y)$ by the mean value theorem, and, since $\|u_\mu\|_\infty\leq C$ by Lemma \ref{Lem:Regularity_u_mu}, we may estimate
	\begin{equation*}
		\begin{split}
			\overline{L_\mu}:=G_\mu(x-y)F(u_\mu(y))&f(u_\mu(x))\varphi(x)\chi_{\{|x-y|>1\}}(x,y)\\
			&\les(|x|+|y|)|u_\mu(y)|^\frac Ns|u_\mu(x)|^{\frac Ns-1}|\varphi(x)|\chi_{\{|x-y|>1\}}(x,y)\\
			&\les\frac{(|x|+|y|)|\varphi(x)|\chi_{\{|x-y|>1\}}(x,y)}{\left(1+|y|^{a\frac Ns}\right)\left(1+|x|^{a\left(\frac Ns-1\right)}\right)}\\
			&\les\frac{\varphi(x)}{1+|y|^{a\left(\frac Ns-1\right)}}\,.
		\end{split}
	\end{equation*}
	Here we used ($f_1$), the decay estimates in Lemma \ref{Lem:Regularity_u_mu} and the compact support of $\varphi$. Since $x\mapsto|x|^{-a\left(\frac Ns-1\right)}$ is integrable in $\R^N\setminus B_1(0)$, thanks to the choice of $a$ in \eqref{w}, by the dominated convergence theorem we conclude that
	\begin{equation}\label{xybig}
		\overline{L_\mu}(x,y)\to\log\frac1{|x-y|}F(u_0(y))f(u_0(x))\varphi(x)\chi_{\{|x-y|>1\}}(x,y)\quad\mbox{in}\ \, L^1(\R^N)\,.
	\end{equation}
\noindent 	Combining \eqref{xysmall} and \eqref{xybig} with \eqref{conv_umu_u0}, we conclude that $u_0$ is a critical point of $J$.
	
\noindent 	The proof of \eqref{logFF} strictly follows the argument in the proof of \cite[Theorem 1.6]{CLR}, being based on Fatou's Lemma, \eqref{improved_integrability_F_mu}, \eqref{unif_bound_norm_u_n_UNIFORMmu}, and Lemma \ref{Lem:MP_level}.
	
\noindent	It remains to show that $u_0$ is nontrivial. Suppose the contrary, then, similarly to \eqref{conv_F} one can show that $\intN f(u_\mu)u_\mu\to0$ and hence
	\begin{equation*}
		\begin{split}
			0&=J_\mu'(u_\mu)[u_\mu]=\|u_\mu\|_V^{\frac Ns}-C_N\intN\left(G_\mu(\cdot)\ast F(u_\mu)\right)f(u_\mu)u_\mu\dd x\\
			&\ges\|u_\mu\|_V^{\frac Ns}-C_N\int\int_{\{|x-y|\leq1\}}\frac{F(u_\mu(y))f(u_\mu(x))u_\mu(x)}{|x-y|^{\tfrac{4(\omega-1)}{3\omega}}}\dd x\dd y\\
			&\geq\|u_\mu\|_V^{\frac Ns}-\intN f(u_\mu)u_\mu=\|u_\mu\|_V^{\frac Ns}+o_\mu(1)\,.
		\end{split}
	\end{equation*}
	This yields $u_\mu\rightharpoonup0$ in $W^{s,\frac Ns}(\R^N)$. But then
	\begin{equation*}
		\begin{split}
			\underline c&\leq c_\mu+o_\mu(1)=J_\mu(u_\mu)\\
			&\leq-C_N\int\int_{\{|x-y|\geq1\}}G_\mu(x-y)F(u_\mu(y))F(u_\mu(x))\dd x\dd y+o_\mu(1)\\
			&\leq C_N\intN\intN\log|x-y|F(u_\mu(y))F(u_\mu(x))\dd x\dd y+o_\mu(1)\\
			&\les\|F(u_\mu)\|_1\intN|x|F(u_\mu(x))\dd x+o_\mu(1)\\
			&\les R\|F(u_\mu)\|_1^2+\int_{\{|x|\geq R\}}\frac{|x|}{1+|x|^{a\frac Ns}}\dd x+o_\mu(1)\leq\varepsilon+o_\mu(1),
		\end{split}
	\end{equation*}
	by \eqref{conv_F} and choosing a sufficiently large $R$. This is a contradiction, and we conclude the proof.
\end{proof}

\section{From Choquard to Schr\"odinger-Poisson: Proof of Theorem \ref{thm_equiv}}\label{Section_SP}
\noindent Theorem \ref{Thm_Chs} yields a positive solution $u$ to the Choquard equation \eqref{Chs}. Next we show that the pair $(u,\phi_u)$, where we recall that $\phi_u:=C_N\log\frac1{|\cdot|}\ast F(u)$, is indeed a solution of the system \eqref{SPs} in the sense of Definition \ref{sol_Poisson}. This step is often neglected in the literature \cite{CW,DW,CW2}. In fact this can be done rigorously as follows, based on a characterisation of the distributional solutions of the Poisson equation by Hyder \cite{H}. The idea is to compare $\phi_u$ with a model solution which we known to solve in the sense of distributions the Poisson equation in \eqref{SPs}, proving that the two solutions may differ only by a constant. The key point is to show that
\begin{equation}\label{logF}
	\intN\log(1+|x|)F(u(x))\dd x\leq C\,,
\end{equation}
which mimics the logarithmic-weighted Poho\v zaev-Trudinger inequality in \cite{CT,BCT}, but just for solutions.

\begin{proof}[Proof of Theorem \ref{thm_equiv}]
	As in the proof of Lemma \ref{Lem:Regularity_u_mu}, following a Nash-Moser iteration technique, it is possible to prove that any positive weak solution of \eqref{Chs} belongs to $L^\infty(\R^N)$ and that there exists $R>0$ such that
	\begin{equation}\label{decay_estimate_u}
		u(x)\les\frac1{1+|x|^\frac{s(2N+3)}{2(N-s)}}\quad\mbox{for all}\ \,|x|>R\,.
	\end{equation}
	Note that the proof of \eqref{decay_estimate_u} differs from the one in Lemma \ref{Lem:Regularity_u_mu} by the kernel, which now is logarithmic; however, since $G_\mu(\cdot)$ is approximating $\log\tfrac1{|\cdot|}$, the strategy to get \eqref{decay_estimate_u} is very similar. From \eqref{decay_estimate_u} and the boundedness in $\R^N$ of $u$, by ($f_1$) one may show that
	\begin{equation}\label{decay_estimate_Fu}
		F(u(x))\les\frac1{1+|x|^{\frac{N(2N+3)}{2(N-s)}}}\quad\mbox{for all}\ \,|x|>R\,.
	\end{equation}
	Therefore
	\begin{equation*}
		\begin{split}
			\intN\log(1+|x|)F(u)\dd x&\leq\|F(u)\|_\infty\int_{\{|x|<R\}}\log(1+|x|)\dd x+C\int_{\{|x|\geq R\}}\frac{\log(1+|x|)}{1+|x|^\frac{N(2N+3)}{2(N-s)}}\dd x\\
			&\leq C+C\int_R^{+\infty}\rho^{N-1-\frac{N(2N+3)}{2(N-s)}}\log\rho\dd\rho<+\infty\,,
		\end{split}
	\end{equation*}
	that is \eqref{logF}. With this in hands, we next show that $\phi_u\in L_\gamma(\R^N)$ for all $\gamma>0$. Indeed,
	\begin{equation*}
\begin{split}
			C_N^{-1}&\intN\frac{|\phi_u(x)|}{1+|x|^{N+2\gamma}}\dd x\leq\intN F(u(y))\left(\intN\left|\log\frac1{|x-y|}\right|\frac1{1+|x|^{N+2\gamma}}\dd x\right)\dd y\\
			&\leq\intN F(u(y))\left(\int_{\{|x-y|>1\}}\frac{\log|x-y|}{1+|x|^{N+2\gamma}}\dd x+\int_{\{|x-y|\leq1\}}\log\frac1{|x-y|}\dd x\right)\dd y\\
			&\leq\intN F(u(y))\bigg(\intN\frac{\log(1+|x|)}{1+|x|^{N+2\gamma}}\dd x\\
			&\quad+\log(1+|y|)\intN\frac{\dd x}{1+|x|^{N+2\gamma}}+\|\log(\cdot)\|_{L^1(B_1(0))}\bigg)\dd y\\
			&\les\intN F(u(y))\left(1+\log(1+|y|)\right)\dd y<+\infty\,,
		\end{split}
	\end{equation*}
	thanks to \eqref{logF} and to the simple inequality $\log|x-y|\leq\log(1+|x|)+\log(1+|y|)$.
	
	\noindent By \cite[Lemma 2.3]{H} the function
	$$\tv_u(x):=C_N\intN\log\left(\frac{1+|y|}{|x-y|}\right)F(u(y))\dd y$$
	belongs to $W^{N-1,1}_{\loc}(\R^N)$ and solves $(-\Delta)^{\frac N2}\tv_u=F(u)$ in $\R^N$ in the sense of Definition \ref{sol_Poisson}. Moreover, by \eqref{logF} it is easy to verify that
	$$\phi_u(x)=\tv_u(x)+C_N\intN\log(1+|y|)F(u(y))\dd x=\tv_u(x)+C\,.$$
	This guarantees that $\phi_u\in W^{N-1,1}_{\loc}(\R^N)$ and solves \eqref{SPs} in the sense of Definition \ref{sol_Poisson} by \cite[Lemma 2.4]{H}, for which all such solutions of $(-\Delta)^{\frac N2}\phi=\mathfrak f$ in $\R^N$ are of the form $\phi=\tv_u+p$ with $p$ polynomial of degree at most $N-1$. The decay behaviour for $\phi_u$ in \eqref{asympt_phi} can be proved following the approach of \cite{CW}, as done in the case $N=2$ in the proof of \cite[Theorem 1.4]{CLR}. It remains to show the H\"older continuity of $u$. To this aim we take advantage of the following local estimate of the H\"older seminorm, obtained in \cite[Corollary 5.5]{IMS}, provided $|(-\Delta)^s_{\frac Ns}u|\leq K$ weakly in $B_{2R}(x_0)$:
	\begin{equation}\label{Holder_estimate_IMS}
		\widetilde CR^\nu[u]_{C^\nu(B_R(x_0))}\leq (KR^N)^{\frac s{N-s}}+\|u\|_{L^\infty(B_{2R}(x_0))}+R^N\left(\int_{\R^N\setminus B_{2R}(x)}\frac{|u(y)|^{\frac Ns-1}}{|x-y|^{2N}}\dd y\right)^{\frac s{N-s}},
	\end{equation}
	where $\widetilde C$ is a universal constant. Note that, since $u\in L^\infty(\R^N)$, it is easy to see that the right-hand side of \eqref{Holder_estimate_IMS} is bounded with a constant depending just on $R,N,s,K$. The assumption $|(-\Delta)^s_{\frac Ns}u|\leq K$ may be verified with similar computations as in \cite[Lemma 3.2]{CLR}, for the bound from above, and, basing on \eqref{logF}-\eqref{decay_estimate_u}, as in \cite[Lemma 3.6]{CLR} for the one from below.
\end{proof}

\end{document}